\RequirePackage{fix-cm}
\pdfoutput=1
\documentclass[leqno, fleqn, centertags, 12pt]{article}

\usepackage{latexsym}
\usepackage{amsmath}
\usepackage{amscd}
\usepackage{amssymb}
\usepackage{verbatim}

\usepackage[T1]{fontenc}
\usepackage{latexsym}
\usepackage{manfnt}
\usepackage[sans]{dsfont}
\usepackage{palatino}
\usepackage[sc, osf]{mathpazo}
\usepackage{eulervm}
\usepackage{euscript}

\usepackage{textcomp}

\usepackage{fullpage}

\makeatletter
\renewcommand{\@makefntext}[1]{\vspace*{0.5ex}\parindent=0em
\hspace*{-0.4em}
\hbox to 0.4em{\hss\@makefnmark}\hspace*{0.4em}{#1}
}
\makeatother

\newcounter{mysectionnumber}
\newcommand{\mysection}[2]{
\setcounter{equation}{0}
\refstepcounter{mysectionnumber}
\section*{ \textnormal{{\themysectionnumber.} {#1}}}\label{#2}}

\numberwithin{equation}{section}

\newcommand{\mynonumbersection}[2]{
\vspace{-0.0ex}
\section*{{}\hspace*{0.00em}$\phantom{1.}$\textnormal{{#1}}}\label{#2}}  

\newcommand{\myit}[1]{\textbf{\textit{#1}}\hspace{0.0em}}

\newcounter{myparnum}[mysectionnumber]
\renewcommand{\themyparnum}{\themysectionnumber.\arabic{myparnum}}
\newcommand{\mypar}[2]{\refstepcounter{myparnum}{\vspace{\medskipamount}\textbf{{\themyparnum. #1}\label{#2}}\hspace{0.5em}}}

\newcommand{\myuppar}[1]{\vspace{\medskipamount}\textbf{#1}\hspace*{0.5em}}

\newcommand{\proof}{\vspace{\medskipamount}{\textbf{{\emph{Proof}.}}\hspace*{0.7em}}}
\newcommand{\eproof}{ $\blacksquare$}

\newcommand{\dis}{\displaystyle}

\def\sss{\hspace{0.05em}\ }
\def\dss{\hspace{0.1em}\ }
\def\trs{\hspace{0.15em}\ }
\def\qss{\hspace{0.2em}\ }
\def\pss{\hspace{0.3em}\ }
\def\oss{\hspace{0.4em}\ }

\def\halfff{\hspace*{0.025em}}
\def\fff{\hspace*{0.05em}}
\def\dff{\hspace*{0.1em}}
\def\trf{\hspace*{0.15em}}
\def\qff{\hspace*{0.2em}}
\def\pff{\hspace*{0.3em}}
\def\off{\hspace*{0.4em}}

\newcommand{\hnsp}{\hspace*{-0.05em}}
\newcommand{\nsp}{\hspace*{-0.1em}}
\newcommand{\nnsp}{\hspace*{-0.15em}}

\newcommand{\dnsp}{\hspace*{-0.2em}}

\newcommand{\id}{\mathop{\mbox{id}}\nolimits}

\newcommand{\cc}{\mathbold{C}}

\newcommand{\rr}{\mathbold{R}}

\newcommand{\kkk}{\mathbold{K}}
\newcommand{\kk}{\mathbold{k}}
\newcommand{\bkk}{\kk^{\dff\perp}}
\newcommand{\nnn}{\mathbold{N}}
\newcommand{\aaa}{\mathbold{A}}

\newcommand{\homv}{\operatorname{Hom}_{\fff a}\fff(V,\qff V)}

\newcommand{\re}{\mbox{Re}\qff}

\newcommand{\toto}{\longrightarrow}
\newcommand{\ttoo}{\hspace*{0.2em}\longrightarrow\hspace*{0.2em}}

\begin{document}

\setlength{\baselineskip}{12pt plus 0pt minus 0pt}
\setlength{\parskip}{12pt plus 0pt minus 0pt}
\setlength{\abovedisplayskip}{12pt plus 0pt minus 0pt}
\setlength{\belowdisplayskip}{12pt plus 0pt minus 0pt}

\newskip\smallskipamount \smallskipamount=3pt plus 0pt minus 0pt
\newskip\medskipamount   \medskipamount  =6pt plus 0pt minus 0pt
\newskip\bigskipamount   \bigskipamount =12pt plus 0pt minus 0pt

\title{{The\qss conformal\qss dilatation\qss
and\qss Beltrami\qss forms\qss\\ 
over\qss quadratic\qss field\qss extensions}}
\date{}
\author{Nikolai V. Ivanov}

\footnotetext{\hspace*{-0.65em}\copyright\ 
Nikolai V. Ivanov, 2017.\trs 
Neither the work reported in this paper,\qss 
nor its preparation were 
supported by any governmental 
or non-governmental agency,\qss 
foundation,\qss 
or institution.}

\maketitle

\myit{\hspace*{0em}\large Contents}\vspace*{1ex} \vspace*{\bigskipamount}\\ 
\myit{Introduction} \hspace*{0.5em} \vspace*{1ex}\\
\myit{\phantom{1}1.}\hspace*{0.5em} Bilinear and quadratic forms \hspace*{0.5em} \vspace*{0.25ex}\\
\myit{\phantom{1}2.}\hspace*{0.5em} The norm and the trace \hspace*{0.5em} \vspace*{0.25ex}\\
\myit{\phantom{1}3.}\hspace*{0.5em} Quadratic forms in dimension $1$ over $\kkk$ \hspace*{0.5em}\vspace*{0.25ex}\\
\myit{\phantom{1}4.}\hspace*{0.5em} The structure of a $\kkk$\dnsp-vector space  on $\aaa\fff(V)$  \hspace*{0.5em} \vspace*{0.25ex}\\
\myit{\phantom{1}5.}\hspace*{0.5em} Conformal structures  \hspace*{0.5em} \vspace*{0.25ex}\\
\myit{\phantom{1}6.}\hspace*{0.5em} Anti-linear maps  \hspace*{0.5em} \vspace*{0.25ex}\\
\myit{\phantom{1}7.}\hspace*{0.5em} The structure of $\homv$  \hspace*{0.5em} \vspace*{0.25ex}\\
\myit{\phantom{1}8.}\hspace*{0.5em} The conformal dilatation and the Beltrami forms \hspace*{0.5em} \vspace*{0.25ex}\\
\myit{\phantom{1}9.}\hspace*{0.5em} Identification of $\mathbb{M}\fff(V)$ with $\homv$   \hspace*{0.5em} \vspace*{0.25ex}\\
\myit{10.}\hspace*{0.5em} Comparing $\mathbold{c}\fff(f)$ and $\mu_{\dff f}$   \hspace*{0.5em} \vspace*{1ex}\\
\myit{References}\hspace*{0.5em}  \hspace*{0.5em}  \vspace*{0.25ex}

\newcommand{\gal}{\operatorname{Gal}\dff}
\newcommand{\tr}{\operatorname{t{\halfff}r}\fff}
\newcommand{\Tr}{\operatorname{T{\halfff}r}\fff}
\newcommand{\oo}[1]{\overline{#1}}

\renewcommand{\baselinestretch}{1}
\selectfont

\mynonumbersection{Introduction}{introduction}

\vspace*{\medskipamount}
\myuppar{The motivation.}
The present paper grew out of the desire to present an invariant treatment of
the geometric interpretation of the complex dilatation suggested in\qss \cite{i-d}\fff.\oss
The main part of\qss \cite{i-d}\qss is devoted to the real linear maps\qss
$\cc\toto \cc$\qss and heavily uses the standard basis\qss $1\fff,\pff i$\qss
of\dss $\cc$\dss as a real vector space.\oss
Only the last section of\qss \cite{i-d}\qss is devoted to
the general case of real linear maps\qss $V\toto W$\dnsp,\oss
where\qss $V\fff,\pff W$\qss are two complex vector spaces of dimension $1$\nnsp,\oss
and this case is treated by a reduction to the case\qss $\cc\toto \cc$\nnsp.\oss
This approach 
served well to the goals
of the paper\qss \cite{i-d}\fff;\oss
see\qss \cite{i-d}\fff,\oss the end of the Introduction.

A recent paper of M.\qss Atiyah\qss \cite{at}\qss reminded the author
that a truly invariant approach requires more than just dealing with maps\qss
$V\toto W$\qss from the very beginning\halfff.\oss
Namely,\oss to quote M.\qss Atiyah,\oss
\emph{``one should not distinguish between\qss $i$\qss and\qss $-\qff i$\qss 
whereas one can distinguish between\qss $1$\qss and\qss $-\qff 1$\nnsp''}.\oss
An attempt to follow this maxim quickly lead to the realization that\pss
\emph{one should distinguish}\pss between\qss $-\qff 1$\qss 
as a purely additive notion and\qss $i^{\fff 2}$\dnsp.\oss

Since\qss $i^{\fff 2}\qff =\qff -\qff 1$\qss anyhow,\oss
the only sensible way to distinguish\qss $-\qff 1$\qss and\qss $i^{\fff 2}$\qss
is to abandon the axiom\qss $i^{\fff 2}\qff +\qff 1\qff =\qff 0$\qss and 
replace\qss $\cc\fff,\pff \rr$\qss by two fields\qss $\kkk\fff,\pff \kk$\qss
such that\dss $\kkk$\dss is an extension of\dss $\kk$\dss of 
degree\dss $[\dff\kkk:\kk\dff]\qff =\qff 2$\dnsp.\oss
Since quadratic forms behave differently in characteristic $2$\nnsp,\oss
it is reasonable to assume that the characteristic of\dss $\kk$\dss is\qss $\neq\qff 2$\dnsp.\oss

This paper is devoted to an analogue of the theory of\qss \cite{i-d}\qss
in this situation,\oss
independent on any choices of bases.\oss
Working with a general field extension\dss $\kkk/\fff\kk$\dss 
does not lead to any new difficulties compared to the classical case\qss
$\kkk/\fff\kk\qff =\qff \cc/\fff\rr$\nnsp,\oss
but only clarifies the theory.\oss 
In the rest of the introduction is devoted to an outline of this analogue.

\myuppar{Conformal structures and conformal dilatation.}
For a vector space $V$ over $\kk$ we denote by\dss $Q\fff(V)$\dss
the vector space of all quadratic forms on $V$\dnsp.\oss
A\qss \emph{conformal structure}\qss on $V$ 
is a non-zero quadratic form on $V$ considered up to 
multiplication by non-zero elements of $\kk$\nnsp,\oss
i.e.\qss a point in
the projective space\qss
$\mathbb{P}\hnsp Q\fff(V)$\qss 
associated with\dss $Q\fff(V)$\dnsp.\oss
If\qss $\dim V\qff =\qff 2$\nnsp,\oss
then\dss $\mathbb{P}\hnsp Q\fff(V)$\dss is a projective plane and
the set\dss $\mathbb{C}\fff(V)\dff\subset\dff \mathbb{P}\hnsp Q\fff(V)$\qss
of conformal structures corresponding to degenerate quadratic forms
is a\qss \emph{conic}\qss in\dss $\mathbb{P}\hnsp Q\fff(V)$\dnsp.\oss
In other terms,\qss $\mathbb{C}\fff(V)$\qss is defined by a homogeneous equation
of degree $2$\nnsp.

This paper is concerned mostly with
vector spaces $V$ of dimension $1$ over $\kkk$\nnsp,\oss
considered as a vector spaces over $\kk$\nnsp.\oss
There is a\qss \emph{canonical conformal structure}\qss 
on any such\dss $V$\dnsp,\oss
denoted by\dss $c_{\fff V}$\nnsp.\oss
It is invariant
under the multiplication by non-zero elements of\dss $\kkk$\nnsp.\oss

Let\dss $\mathbb{A}\fff(V)$\dss be the projective line polar\dff\footnotemark\qss
to\dss $c_{\fff V}$\dss 
with respect to the conic\dss $\mathbb{C}\fff(V)$\dnsp.\oss
Then 
\[
\quad
\mathbb{M}\fff(V)
\off =\off
\mathbb{P}\hnsp Q\fff(V)\qff \smallsetminus\qff \mathbb{A}\fff(V)
\]
is an affine plane containing\dss $c_{\fff V}$\nnsp.\oss
One can turn\dss $\mathbb{M}\fff(V)$\dss into a vector space over $\kk$
by designating\dss $c_{\fff V}$\dss as its zero vector\halfff.\oss
The structure of a vector space of dimension\dss $1$\dss over\dss $\kkk$\dss on\dss $V$\dss
allows to extend this structure of a vector space over\dss $\kk$\dss 
on\dss $\mathbb{M}\fff(V)$\dss to a structure of a vector space over\dss $\kkk$\dss
on\dss $\mathbb{M}\fff(V)$\dss
in a canonical manner\halfff.\oss
There is a canonical quadratic form\oss 
$\dis
\mathbb{D}\dff\colon\dff \mathbb{M}\fff(V)\toto \kk
$\oss
equal to\dss 
$1$\dss on the conic\qss $\mathbb{C}\fff(V)$\dnsp.\oss

\footnotetext{In the main part of the paper the notion of polarity is not used 
and\dss $\mathbb{A}\fff(V)$\dss 
is defined differently\halfff.\oss
Theorem\qss \ref{description-of-complement}\qss implies that
the two definitions are equivalent\halfff.}

Let\qss $f\dff\colon\dff V\toto W$\qss be a map linear over\dss $\kk$\dss between vector spaces\qss
$V\fff,\pff W$\qss of dimension\dss $1$\dss over\dss $\kkk$\nnsp.\oss
A natural measure of the distortion of the canonical conformal structure by\dss $f$\dss is
the pull-back\oss
$\dis
\mathbold{c}\fff(f)
\off =\off
f^{\dff *}\fff c_{\fff W}
\qff \in\qff
\mathbb{P}\hnsp Q\fff(V)$\dnsp,\qff\oss
called\dss the\qss \emph{conformal dilatation}\qss of\qss $f$\nnsp.\oss

\myuppar{Anti-linear maps and Beltrami forms.}
Since\dss \dss $[\dff\kkk:\kk\dff]\qff =\qff 2$\nnsp,\oss
there is only one non-trivial automorphism of\dss $\kkk$ fixed on $\kk$\nnsp,\oss
which is denoted by\qss $z\qff \longmapsto\qff \oo{z}$\nnsp.\oss

Let\qss $V\fff,\pff W$\qss be vector spaces of dimension $1$ over $\kkk$\nnsp.\oss
A map\dss $g\dff\colon\dff V\toto W$\dss 
is said to be\qss \emph{anti-linear}\qss if
it is linear over $\kk$ and\qss
$\dis
g\fff(z\dff v)
\pff =\pff 
\oo{z}\qff g\fff(v)$\qss
for all\qss $z\dff\in\dff \kkk$\nnsp,\qss $v\dff\in\dff V$\dnsp.\oss
Such maps form a vector space\qss 
$\operatorname{Hom}_{\fff a}\fff(V\fff,\trf W)$\qss  over $\kkk$\nnsp.\oss
There is a canonical quadratic form\qss 
\[
\quad
\mathcal{D}\dff\colon\dff\homv\toto \kk\dff.
\]
In view of Lemma\qss \ref{det-trace}\qss it can be defined by\oss 
$\mathcal{D}\fff(g)\qff =\qff -\qff \det g$\nnsp,\oss
where the determinant is taken over $\kk$\nnsp.\oss

If a map\qss $f\dff\colon\dff V\toto W$\qss is linear over $\kk$\nnsp,\oss
then\oss 
$\dis
f
\off =\off
Lf\qff +\qff Af$\nnsp,\oss
where\dss $Lf$\dss is linear over\dss $\kkk$\dss
and\dss $Af$\dss is anti-linear\halfff.\oss
Moreover\halfff,\oss the maps\qss $Lf\fff,\pff Af$\qss
are uniquely determined by $f$\nnsp.\oss
If\qss $Lf\qff \neq\qff 0$\dnsp,\oss
then\dss $Lf$\dss is invertible and the composition\qss
\[
\quad
\mu_{\dff f}
\off =\off
(Lf)^{\fff -\dff 1}\circ (Af)
\qff \in\qff \homv\dff.
\]
is a natural measure of failure of $f$ to
be linear over $\kkk$\nnsp,\oss
called the\qss \emph{Beltrami form}\dss of\dss $f$\nnsp.\oss

\myuppar{Comparing the conformal dilatation and Beltrami forms.}
The results of\qss \cite{i-d}\qss suggest that\dss $\mathbold{c}\fff(f)$\dss
and\dss $\mu_{\dff f}$\dss should be related by an analogue of the map
relating the Klein and the Poincar\'{e} models of the hyperbolic plane.\oss
As it turns out\halfff,\oss this is indeed the case.

There is a canonical isomorphism of\dss $\mathbb{M}\fff(V)$\dss and\dss $\homv$\dss
as vector spaces over $\kkk$\nnsp.\oss
Moreover\halfff,\oss this isomorphism turns
the quadratic forms\dss $\mathbb{D}$\nnsp,\qss $\mathcal{D}$\dss 
one into the other\halfff.\oss
One can use this isomorphism to identify\dss $\mathbb{M}\fff(V)$\dss and\dss $\homv$\dnsp.\oss

Let\qss 
$f\dff\colon\dff V\toto W$\qss be a map linear over\qss $\kk$\nnsp.\oss 
The map\dss $f$\dss is called\qss \emph{regular}\pss
if\qss $\mathbold{c}\fff(f)\dff\in\dff \mathbb{M}\fff(V)$\dnsp,\oss
and\qss \emph{exceptional}\qss otherwise.\oss
The exceptional maps are a new feature of the general case.\oss
If the extension\dss $\kkk/\fff\halfff\kk$\dss resembles enough
the classical case of\dss $\cc/\fff\rr$\nnsp,\oss
then all non-zero maps\dss $f$\dss 
are regular\halfff.\oss
This is the case if the field $\kk$ is ordered and $\kkk$
is obtained by adjoining a root of a\qss \emph{negative}\qss element to $\kk$\nnsp,\oss
as it easily follows from Theorem\qss \ref{main-regular}.\oss

Suppose that\qss $Lf\qff \neq\qff 0$.\qff\oss 
If\oss 
$\dis
1\qff +\qff \mathcal{D}\fff(\mu_{\dff f})\off \neq\off 0$\dnsp,\oss
then\dss $f$\dss is regular and 
\[
\quad
\mathbold{c}\fff(f)
\off\qff =\off\qff
\frac{2\dff \mu_{\dff f}}{\fff 1\qff +\qff \mathcal{D}\fff(\mu_{\dff f})\fff} 
\]
after the identification of\pss $\mathbb{M}\fff(V)$\qss
with\qss $\homv$\dnsp.\qff\oss
See Theorem\qss \ref{main-regular}.\oss
A simple calculation\qss 
(see\oss \cite{i-d}\fff,\oss Corollary\qss 3.3)\qss 
shows that for\qss 
$\kkk\qff =\qff \cc$\nnsp,\qss $\kk\qff =\qff \rr$\nnsp,\oss 
and\qss $V\qff =\qff \cc$\qss
this relation between\dss $\mathbold{c}\fff(f)$\dss and\dss $\mu_{\dff f}$\qss
is the same as in\oss \cite{i-d}.\qff\oss

If\oss
$\dis
1\qff +\qff \mathcal{D}\fff(\mu_{\dff f})\off =\off 0$\dnsp,\qff\oss
then\dss $f$\dss is exceptional.\oss
In this case the relation between\dss $\mathbold{c}\fff(f)$\dss
and\dss $\mu_{\dff f}$\dss is even simpler\halfff.\oss
See Theorem\qss \ref{main-exceptional}.\oss

\mysection{Bilinear\qss and\qss quadratic\qss forms}{quadratic forms}

\vspace*{\medskipamount}
\myuppar{Bilinear forms,\qss their matrices and determinants.}
Let us fix a vector space $V$ over $\kk$ of finite dimension $n$
and a basis\qss $v_1\fff,\pff v_2\fff,\pff \ldots\fff,\pff v_{\fff n}$\qss of $V$\dnsp.\oss
A bilinear form\qss $\varphi\colon V\times V\toto \kk$\qss gives rise to 
an\qss $n\times n$\qss matrix\qss 
$M\fff(\varphi)\qff =\qff (\fff M_{\dff i\dff j}\fff(\varphi)\fff)$\trs with the entries
\[
\quad
M_{\dff i\dff j}\fff(\varphi)
\off =\off
\varphi\fff(\fff v_{\fff i}\dff,\pff v_{\fff j}\fff)\dff,
\]
called\qss \emph{the matrix of}\dss $\varphi$\dss with respect to the basis\qss
$v_1\fff,\pff v_2\fff,\pff \ldots\fff,\pff v_n$\nnsp.\oss
The determinant
\[
\quad
\det \varphi \off =\off \det M\fff(\varphi)
\]
is called the\qss \emph{determinant}\qss of\dss $\varphi$\dss
(with respect to the basis\qss $v_1\fff,\pff v_2\fff,\pff \ldots\fff,\pff v_{\fff n}$\nsp).\oss

\myuppar{Change of basis.}
If\qss $v'_1\fff,\pff v'_2\fff,\pff \ldots\fff,\pff v'_{\fff n}$\qss 
is some other basis of $V$\dnsp,\oss then
\[
\quad
v'_{\dff i}
\off =\off
a_{\dff i\dff 1}\dff v_1\qff +\qff a_{\dff i\dff 2}\dff v_2\qff +\qff 
\ldots\qff +\qff a_{\dff i\dff n}\dff v_{\fff n}
\]
for some invertible matrix\qss 
$A\qff =\qff (a_{\dff i\dff j}\fff)$\qss
and\oss
\[
\quad
M'\fff(\varphi)
\off =\off
A\dff M\fff(\varphi)\dff A^{t}\dff,
\]
where\dss $M'\fff(\varphi)$\dss is the matrix of\dss $\varphi$\dss 
with respect to the basis\qss
$v'_1\fff,\pff v'_2\fff,\pff \ldots\fff,\pff v'_{\fff n}$\nnsp,\oss
and\qss $A^{t}$\qss is the transposed matrix of\dss $A$\nnsp.\oss
Therefore\oss
\[
\quad
{\det}{\fff '} \varphi\off =\off \det \varphi \cdot (\det A)^{\fff 2}\dff,
\]
where\dss $\det' \varphi$\dss is the determinant of\dss $\varphi$\dss
with respect to the basis\qss
$v'_1\fff,\pff v'_2\fff,\pff \ldots\fff,\pff v'_{\fff n}$\nnsp.\oss

\myuppar{Bilinear forms and linear maps.}
If\qss $f\colon V\toto V$\qss is a linear map,\oss
then\oss
\[
\quad
f\fff(v_{\dff i})
\off =\off
f_{\dff i\dff 1}\dff v_1\qff +\qff f_{\dff i\dff 2}\dff v_2\qff +\qff 
\ldots\qff +\qff f_{\dff i\dff n}\dff v_{\fff n}\dff,
\]
where\qss $F\qff =\qff (f_{\dff i\dff j}\fff)$\qss is the matrix of\dss $f$\dss 
in the basis\qss
$v_1\fff,\pff v_2\fff,\pff \ldots\fff,\pff v_{\fff n}$\nnsp.\oss
Together with the form\dss $\varphi$\dss the map\dss $f$\dss gives rise to a new bilinear form,\oss
the\qss \emph{left\dss product\qss
\[
\quad
f\nsp\cdot\nsp\varphi\qff \colon\qff
(v\fff,\pff w)\qff \longmapsto\qff \varphi\fff(\fff f\fff(v)\fff,\pff w\fff)\dff,
\]
of\qss $f$\trs and\trs $\varphi$}\nnsp.\oss
The matrix\dss
$M\fff(f\nsp\cdot\nsp\varphi)$\dss of\dss $f\nsp\cdot\nsp\varphi$\dss is equal to\dss 
$F\dff M\fff(\varphi)$\dnsp,\oss and hence\oss
\begin{equation}
\label{dot-det}
\quad
\det\qff (f\nsp\cdot\nsp\varphi)
\off =\off
\det\dff f\cdot \det\dff \varphi\dff,
\end{equation}
where\qss $\det f\qff =\qff \det F$\qss is independent on the choice of the basis.

\myuppar{Non-degenerate bilinear forms.}
A bilinear form\dss $\varphi$\dss gives rise to a linear map
\[
\quad
\widehat{\varphi}\qff \colon\qff
V\ttoo \hom (V\fff,\pff \kk)
\]
taking\qss $v\dff\in\dff V$\qss to the map\qss
$w\qff\longmapsto\qff \varphi\fff(v\fff,\pff w)$\dnsp.\oss
The form\dss $\varphi$\dss is said to be\qss \emph{non-degenerate}\pss
if\dss $\widehat{\varphi}$\dss is an isomorphism.\oss
As is well known,\qss this condition is equivalent to\qss
$\det \varphi\qff \neq\qff 0$\dnsp.\oss

Suppose that\dss $\varphi\fff,\pff \beta$\dss are bilinear forms on $V$ and\dss
$\varphi$\dss is non-degenerate.\oss
By the previous paragraph,\oss
for every\dss $v\dff\in\dff V$\dss the homomorphism\qss
$w\qff\longmapsto\qff \beta\fff(v\fff,\pff w)$\qss
is equal to\qss
\[
\quad
w\qff\longmapsto\qff \varphi\fff(f\fff(v)\fff,\pff w)
\]
for a unique\qss $f\fff(v)\dff\in\dff V$\dnsp.\oss
A trivial verification shows that\dss $f$\dss is linear
and hence\qss $\beta\qff =\qff f\nsp\cdot\nsp\varphi$\nnsp.

\myuppar{Quadratic forms.}
A bilinear form\qss
$\varphi\colon V\times V\toto \kk$\qss is called\qss \emph{symmetric}\qss if\oss
\[
\quad
\varphi\fff(v\fff,\pff w)
\off =\off 
\varphi\fff(w\fff,\pff v)
\]
for all\qss $v\fff,\pff w\dff\in\dff V$\dnsp.\oss
A bilinear form  
is symmetric if and only if its matrix is symmetric.\oss
A\dss \emph{quadratic form}\qss on $V$ 
is a map\qss $q\colon V\toto \kk$\qss of the form\qss
$v\qff \longmapsto\qff \varphi\fff(v\fff,\pff v)$\dnsp,\oss
where\dss $\varphi$\dss is a 
symmetric bilinear form.\oss 
The form\dss $\varphi$\dss is uniquely determined by $q$ as
its\qss \emph{polarization}\qss 
\[
\quad
(\fff v\fff,\qff w\fff)_q
\off =\off 
\bigl(q\fff(v\qff +\qff w)\qff -\qff q\fff(v)\qff -\qff q\fff(w)\dff\bigr)\bigl/2\qff.
\]
A quadratic form $q$ is said to be\qss \emph{non-degenerate}\qss
if its polarization 
is non-degenerate.\oss

The\qss \emph{Gram matrix}\sss $G\fff(q)$\dss of a quadratic form\dss $q$\dss
is defined as the matrix of its polarization,\oss
and the determinant\dss $\det q$\dss of\dss $q$\dss is defined as the determinant
of its polarization,\oss i.e.
\[
\quad
\det q\off =\off \det G\fff(q)\dff.
\]
A quadratic form\dss $q$\dss is non-degenerate if and only if\qss
$\det q\qff \neq\qff 0$\dnsp.\oss

\myuppar{The determinant map.}
The set\dss $Q\fff(V)$\dss of quadratic forms on $V$ is a vector space over $\kk$\nnsp.\oss
The map\qss $q\qff\longmapsto\qff G\fff(q)$\qss
is an isomorphism between\dss $Q\fff(V)$\dss
and the space of symmetric\qss $n\times n$\qss matrices
with entries in\dss $\kk$\nnsp.\oss
The determinants\dss $\det q$\dss 
lead to the\qss \emph{determinant map}\qss
\[
\quad
\det \colon Q\fff(V)\ttoo \kk\dff.
\]
If\dss $\det\halfff'$\dss is the determinant map corresponding to 
some other basis of\dss $V$,\oss then
\[
\quad
\det\halfff'\fff q\off =\off \det q \cdot (\det A)^{\fff 2}\dff,
\]
where\dss $A$\dss is the matrix defined as above.\oss
It follows that up to a multiplication by a non-zero constant the determinant map
is independent from the choice of basis of\dss $V$\dnsp.\oss

\myuppar{Symmetric maps.}
A linear map\qss
$f\dff \colon\dff V\toto V$\qss 
is said to be\qss \emph{symmetric}\qss with respect to a bilinear form\dss $\varphi$\trs if
\[
\quad
\varphi\fff(\fff f\fff(v)\fff,\pff w\fff)
\off =\off
\varphi\fff(\fff v\fff,\pff f\fff(w)\fff)
\]
for all\qss $v\fff,\pff w\dff\in\dff V$\dnsp.\qff\oss
If\dss $\varphi$\dss is symmetric,\oss then\oss
$\dis
\varphi\fff(\fff v\fff,\pff f\fff(w)\fff)
\pff =\pff
\varphi\fff(\fff f\fff(w)\fff,\pff v\fff)$\oss
for all\qss $v\fff,\pff w\dff\in\dff V$\dnsp.\oss
It follows that\dss $f$\dss is symmetric with respect to a symmetric bilinear form\dss $\varphi$\dss
if and only if the left product\dss $f\nsp\cdot\nsp\varphi$\dss is symmetric.\oss

\myuppar{Quadratic forms and symmetric maps.}
A linear map\qss
$f\dff \colon\dff V\toto V$\qss  
is said to be\qss \emph{symmetric}\qss 
with respect to a quadratic form\dss
$q$\dss on\dss $V$\dss if\dss $f$\dss 
is symmetric with respect to the polarization of\dss $q$\nnsp,\oss
i.e.\qss if\dss the bilinear form
\[
\quad
(\fff v\fff,\pff w\fff)
\qff \longmapsto\qff
(\fff f\fff(v)\fff,\pff w\fff)_q
\]
is symmetric.\oss 
If\dss $q$\dss is a quadratic form with the polarization\dss $\varphi$\dss
and\dss $f$\dss is a linear map symmetric with respect to\dss $q$\nnsp,\oss
then the\qss \emph{product}\qss $f\nsp\cdot\nsp q$\qss  
is defined as the quadratic form with the polarization\dss  $f\nsp\cdot\nsp \varphi$\nnsp.\oss
Equivalently\halfff,\oss the\qss \emph{product of\qss $f$\dss and\dss $q$}\qss is the map
\[
\quad
f\nsp\cdot\nsp q\dff\colon\dff v\qff\longmapsto\qff (\fff f\fff(v)\fff,\qff v\fff)_q\dff.
\]
By applying\qss (\ref{dot-det})\qss to\dss $f$\dss 
and the polarization\dss $\varphi$\dss of\dss $q$\nnsp,\oss
we see that\oss
\begin{equation}
\label{det-dot-quadratic}
\quad
\det\qff (f\nsp\cdot\nsp q)
\off =\off
\det\dff f \cdot \det\dff q\dff.
\end{equation}

\vspace*{-\medskipamount}
Suppose that\dss $n$\dss is a non-degenerate quadratic form on\dss $V$\dnsp.\oss
Let\dss $p$\dss be an arbitrary quadratic form on\dss $V$\dnsp.\oss 
Then there exist a unique linear map\qss $f\dff\colon\dff V\toto V$\qss such that
\[
\quad
(\fff v\fff,\pff w\fff)_{\fff p}
\off =\off
(\fff f\fff(v)\fff,\pff w\fff)_n
\]
for all\qss $v\fff,\pff w\dff\in\dff V$\dnsp.\oss
Since the polarization\qss $(\fff v\fff,\pff w\fff)_{\fff p}$\qss is symmetric,\qss
$f$\dss is symmetric with respect to\dss $n$\nnsp.\oss
By taking\qss $v\qff =\qff w$\qss in the last displayed formula,\oss
we see that\qss $p\qff =\qff f\nsp\cdot\nsp n$\nnsp.\oss
The map\dss $f$\dss is uniquely determined by the forms\qss $p\fff,\pff n$\nnsp.\oss
We will denote it by\dss $p/n$\nnsp.\oss

The operations\qss $p\dff \longmapsto p/n$\qss and 
$f\dff \longmapsto f\nsp\cdot\nsp n$ are mutually inverse in the sense that\oss
\[
\quad
p\off =\off (p/n)\nsp\cdot\nsp n
\hspace*{1.5em}\mbox{ and }\hspace*{1.5em}
f\off =\off f\nsp\cdot\nsp n\fff/n
\]
for all\qss $p\dff\in\dff Q\fff(V)$\qss
and all maps\qss $f\fff\colon\fff V\toto V$\qss
symmetric with respect to\dss $n$\nnsp.\oss

\myuppar{Anisotropic forms.}
A quadratic form\dss $q$\dss is called\qss
\emph{anisotropic}\qss if\qss $q\fff(v)\qff \neq\qff 0$\qss 
for\qss $v\qff \neq\qff 0$\dnsp.\oss
Suppose that\dss $q$\dss is an anisotropic form with the polarization\dss $\varphi$\nnsp.\oss
If\qss $v\qff \neq\qff 0$\dnsp,\oss then the map\qss 
$w\qff\longmapsto\qff \varphi\fff(v\fff,\pff w)$\qss
is non-zero because\qss
$\varphi\fff(v\fff,\pff v)\qff =\qff q\fff(v)\qff \neq\qff 0$\dnsp.\oss
It follows that\dss $\widehat{\varphi}$\dss is injective
and hence is an isomorphism.\oss
Therefore\qss $\varphi$\dss and\dss $q$\dss are non-degenerate.\oss

\myuppar{The case of vector spaces of dimension $2$\nnsp.}
Suppose now that the dimension of\dss $V$\dss is $2$\nnsp.\oss
Let\qss $v\fff,\pff w$\qss be a basis of $V$\dnsp.\qff\oss
If\pss $q\dff\in\dff Q\fff(V)$\dnsp,\qff\oss then
\vspace*{\smallskipamount}
\[
\quad
G\fff(q)
\off =\off
\begin{bmatrix}
\pff(\fff v\fff,\qff v\fff)_q & (\fff v\fff,\qff w\fff)_q\fff\off\vspace*{1.5\medskipamount} \\
\off\qff(\fff w\fff,\qff v\fff)_q & \pff\fff(\fff w\fff,\qff w\fff)_q\fff\off
\end{bmatrix}
\]

\vspace{-0.75\bigskipamount}
is the Gram matrix of $q$ with respect to the basis\qss $v\fff,\pff w$\nnsp.\oss 
It follows that the determinant map\oss
$\dis 
\det \colon Q\fff(V)\ttoo \kk$\oss
with respect to the basis\qss $v\fff,\pff w$\qss has the form
\vspace*{\smallskipamount}
\[
\quad
\det q
\off =\off
(\fff v\fff,\qff v\fff)_q\dff (\fff w\fff,\qff w\fff)_q
\qff -\qff\fff
(\fff v\fff,\qff w\fff)_q^{\fff 2}
\off =\off
q\fff(v)\dff q\fff (w)\qff -\qff (\fff v\fff,\qff w\fff)_q^{\fff 2}\qff.
\]

\vspace*{-0.75\bigskipamount}
Therefore the determinant map itself is a non-degenerate quadratic form.\oss

\myuppar{Orthogonality\halfff.}
Two vectors\dss $v\fff,\pff w\dff\in\dff V$\dss are said to be\qss
\emph{orthogonal}\qss with respect to a form\dss $q\dff\in\dff Q\fff(V)$\dss if\qss
$(\fff v\fff,\qff w\fff)_q\qff =\qff 0$\dnsp,\oss i.e.\qss if
\[
\quad
q\fff(v\qff +\qff w)\off =\off q\fff(v)\qff +\qff q\fff(w)\dff.
\]
Suppose now that the dimension of\dss $V$\dss is\dss $2$\nnsp.\oss
Let\qss $p\fff,\pff q$\qss be two quadratic forms on $V$\dnsp.\oss
They are said to be\qss \emph{orthogonal}\qss
if\qss $p\fff,\pff q$\qss
are orthogonal with respect to 
a determinant map $\det$ considered as
a quadratic form on $Q\fff(V)$\dnsp,\oss i.e.\qss if 
\begin{equation}
\label{orthogonal}
\quad
\det\qff (p\qff +\qff q)
\off =\off
\det p\qff +\qff \det q\dff.
\end{equation}
Since different determinant maps differ by the multiplication
by a non-zero element of $\kk$\nnsp,\oss
the property\qss (\ref{orthogonal})\qss does not depend on the choice
of the determinant map.\oss

\mypar{Theorem.}{det-orthogonality} 
\emph{Let\qss $V$\qss be a vector space over\qss $\kk$\qss of dimension\qss $2$\qss
and\dss let\pss $p\fff,\pff q\dff\in\dff Q\fff(V)$\dnsp.
Suppose that\qss $v\fff,\pff w$\qss is a basis of\qss $V$\qss 
orthogonal with respect to\qss $p$\dnsp.\oss
Then the forms\qss $p\fff,\pff q$\qss are orthogonal with respect
to any determinant map if and only if}
\[
\quad
p\fff(v)\dff q\fff(w)\qff +\qff q\fff(v)\trf p\fff(w)
\off =\off 0\dff.
\]

\vspace*{-0.75\bigskipamount}
\proof\qss By the remark preceding the theorem,\oss
it is sufficient to consider the orthogonality with respect to determinant map 
defined by the basis\qss 
$v\fff,\pff w$\dnsp.\oss

The forms\qss $p\fff,\pff q$\qss are orthogonal
if and only if\qss (\ref{orthogonal})\qss holds,\oss
i.e.\qss if\dss and only if\vspace*{2\smallskipamount}
\begin{equation}
\label{p+q}
\quad
(p\qff +\qff q)\fff(v)\dff (p\qff +\qff q)\fff (w)
\qff -\qff 
(\fff v\fff,\qff w\fff)_{p\qff +\qff q}^{\fff 2}
\off =
\end{equation}

\vspace*{-2.75\bigskipamount}
\[
\quad
\hspace*{8em}
=\off
\Bigl(p\fff(v)\dff p\fff (w)\qff -\qff (\fff v\fff,\qff w\fff)_p^{\fff 2}\Bigr)
\off +\off
\Bigl(q\fff(v)\dff q\fff (w)\qff -\qff (\fff v\fff,\qff w\fff)_q^{\fff 2}\Bigr).
\]

\vspace{-0.75\bigskipamount}
The basis\qss $v\fff,\pff w$\qss is orthogonal with respect to $p$ and hence\qss
$(\fff v\fff,\qff w\fff)_p\qff =\qff 0$\qss and\oss
\[
\quad
(\fff v\fff,\qff w\fff)_{p\qff +\qff q}
\off =\off 
(\fff v\fff,\qff w\fff)_p\qff +\qff (\fff v\fff,\qff w\fff)_q
\off =\off
(\fff v\fff,\qff w\fff)_q\qff.
\]
It follows that the\qss (\ref{p+q})\qss is equivalent to\vspace*{\smallskipamount}
\[
\quad
(p\qff +\qff q)\fff(v)\dff (p\qff +\qff q)\fff (w)
\qff -\qff
(\fff v\fff,\qff w\fff)_q^{\fff 2}
\off =\off
p\fff(v)\dff p\fff(w)
\qff +\qff
q\fff(v)\dff q\fff(w)\qff -\qff (\fff v\fff,\qff w\fff)_q^{\fff 2}\qff,
\]

\vspace{-0.75\bigskipamount}
or\halfff,\oss what is the same,\oss to\vspace*{\smallskipamount}
\begin{equation}
\label{p+q-0}
\quad
(p\qff +\qff q)\fff(v)\dff (p\qff +\qff q)\fff (w)
\off =\off
p\fff(v)\dff p\fff(w)
\qff +\qff
q\fff(v)\dff q\fff(w)\qff.
\end{equation}

\vspace*{-0.75\bigskipamount}
Since\qss $(p\qff +\qff q)\fff(u)\qff =\qff p\fff(u)\qff +\qff q\fff(u)$\qss
for all\qss $u\dff\in\dff V$\dnsp,\pss 
the left hand side of\qss (\ref{p+q-0})\qss is equal to\vspace*{\smallskipamount}
\[
\quad
p\fff(v)\trf p\fff(w)\qff +\qff q\fff(v)\dff q\fff(w)
\qff +\qff
p\fff(v)\dff q\fff(w)\qff +\qff q\fff(v)\trf p\fff(w)
\]

\vspace*{-0.75\bigskipamount}
and hence\qss (\ref{p+q-0})\qss is equivalent to\qss
$\dis
p\fff(v)\dff q\fff(w)\qff +\qff q\fff(v)\trf p\fff(w)
\qff =\qff 0$\dnsp.\pss
The theorem follows.\pss  \eproof

\mysection{The\qss norm\qss and\qss the\qss trace}{norm-trace}

\vspace*{\medskipamount}
\myuppar{The norm and the trace.}
Recall that\qss $z\qff \longmapsto\qff \oo{z}$\qss is the only
non-trivial automorphism of the field $\kkk$ fixed on $\kk$\nnsp.\oss
Let\qss $N\fff,\pff \tr\dff\colon\dff \kkk\toto \kk$\qss 
be the maps
\[
\quad
N\fff(z)\off =\off z\dff\oo{z}\qff,
\hspace*{1.5em}
\tr\fff(z)\off =\off \left(z\qff +\qff \oo{z}\fff\right)\left/2\dff.\right.
\]
Then\dss $N$\dss is the norm of the extension\dss $\kkk/\fff\kk$\dss and\dss
$\tr$\dss is the half of the trace of\dss $\kkk/\fff\kk$\nnsp.\oss
The map\dss 
$\tr$\dss has the advantage of being equal to the identity on\dss $\kk$\nnsp.\oss
If\oss $\kkk/\fff\kk\off =\off \cc/\fff\rr$\nnsp,\oss
then\qss $\oo{z}$\qss is the complex conjugate of\qss $z$\nnsp,\oss
$\tr\fff(z)\off =\off \re z$\oss and\oss
$N\fff(z)\off =\dff\off |\dff z\dff|^{\fff 2}$\nnsp.

\myuppar{The norm as a quadratic form on\dss $\kkk$\nnsp.}\oss
Let\oss
\[
\langle\fff\halfff z\fff,\qff u\fff\rangle\off =\off\dff \tr\dff(z\trf\oo{u}\fff)\dff.
\]
Then\qss 
$\langle\fff\halfff z\fff,\qff u\fff\rangle$\qss 
is a bilinear form on\dss $\kkk$\dss
considered as a vector space over\dss $\kk$\nnsp.\oss
Since
\[
\quad
\tr\dff(z\qff\oo{u}\fff)
\off =\off
\tr\dff(\dff\oo{\fff z\qff \oo{u}\fff}\dff)
\off =\off
\tr\dff(u\qff\oo{z}\fff)
\]
for all\qss $z\fff,\pff u\dff\in\dff \kkk$\nnsp,\oss
the bilinear form\qss
$\langle\fff\halfff z\fff,\qff u\fff\rangle$\qss
is symmetric.\oss
Since\oss
$
z\dff\oo{z}
\off =\off
\tr(z\dff\oo{z})
$\oss
for all\qss $z\dff\in\dff \kkk$\nnsp,\oss
it follows that\dss $N$\dss
is a quadratic form on\dss $\kkk$\dss 
with polarization\qss
$\langle\fff\halfff z\fff,\qff u\fff\rangle$\nnsp.\oss
Obviously,\oss the quadratic form\dss $N$\dss is anisotropic and 
hence is non-degenerate.\oss

\myuppar{The orthogonal complement of\dss $\kk$\nnsp.}
Let\dss $\bkk$\dss be the orthogonal complement of\dss $\kk$\dss
in\dss $\kkk$\dss with respect to\dss $N$\nnsp,\oss
i.e.\qss the set of all\qss $\rho\dff\in\dff \kkk$\qss such that\qss
$\langle\fff \rho\fff,\qff a\fff\rangle
\qff =\qff 
0$\qss for all\qss
$a\dff\in\dff \kk$\nnsp.\oss
Since\dss $N$\dss is anisotropic,\qss
$\kk\dff \cap\dff \bkk\qff =\qff \varnothing$\qss
and hence\qss
$\kkk\qff =\qff \kk\oplus\bkk$\nnsp.\qff\oss
Since the dimension of\dss $\kkk$\dss over\dss $\kk$\dss is equal to\dss $2$\nnsp,\oss
the dimension of\dss $\bkk$\dss over\dss $\kk$\dss is equal to $1$\nnsp.\oss

The orthogonal complement\qss
$\bkk$\dss
is equal to 
the kernel of\qss $\tr\dff\colon\dff\kkk\toto\kk$\nnsp,\oss
as it follows from the fact that\qss 
$\langle\fff \rho\fff,\qff a\fff\rangle
\qff =\qff
a\dff\tr\fff(\rho)$\qss for\qss $a\dff\in\dff \kk$\nnsp.\oss
Therefore\qss
$\rho\dff\in\dff \bkk$\qss
implies that\oss
\[
\quad
\oo{\rho}\off =\qff -\off \rho
\hspace*{1.5em}\mbox{ and }\hspace*{1.5em}
\rho^{\fff 2}
\off =\off 
-\qff \oo{\rho}\dff\rho
\off =\off 
-\qff N\fff(\rho)\dff\in\dff \kk\dff.
\]

\vspace*{-\bigskipamount}
\mysection{Quadratic\qss forms\qss in\qss dimension\qss $1$\qss over\qss $\kkk$}{forms-over-extension}

\vspace*{\medskipamount}
\myuppar{The framework.} \emph{For the remaining part of the paper we will
assume that\dss $V$\dss is a vector space over\dss $\kkk$\dss of dimension\dss $1$\dss over\dss $\kkk$\nnsp.}\oss
By the restriction of scalars we may consider $V$ also as a vector space over $\kk$\nnsp.\oss
By a quadratic form on $V$ we will understand a quadratic form on $V$
as a vector space over $\kk$\nnsp.\oss
As a vector space over $\kk$
the vector space $Q\fff(V)$ of quadratic forms on $V$
depends only on the structure of $V$ as a vector space over $\kk$\nnsp.\oss

\myuppar{The norm-like quadratic forms.}
A form\qss $n\dff\in\dff Q\fff(V)$\qss is said to be\qss \emph{norm-like}\qss if
\begin{equation}
\label{n-forms}
\quad
n\fff(z\dff v)
\off =\off
N\fff(z)\trf n\fff(v)
\end{equation}
for all\qss $v\dff\in\dff V$\qss and\qss $z\dff\in\dff \kkk$\nnsp.\oss
Let\dss $\nnn\fff(V)$\dss be the space of all norm-like forms on\dss $V$\dnsp.\oss

The dimension of\qss $\nnn(V)$\qss over\qss $\kk$\qss is equal to\qss $1$\nnsp.\oss
Indeed,\oss let\qss $v\dff\in\dff V$\qss and\qss $v\qff \neq\qff 0$\dnsp.\oss
Then every element of\dss $V$\dss is equal to\dss $z\dff v$\dss for some\qss
$z\dff\in\dff \kkk$\qss and hence every norm-like form is determined by its value at\dss $v$\nnsp.\oss
On the other hand,\oss
for every\dss $a\dff\in\dff \kk$\dss the map\oss
$\dis
z\dff v\dff \longmapsto\dff N\fff(z)\dff a$\oss 
is a norm-like form and takes the value $a$ at $v$\nnsp.\oss

Every non-zero\dss $n\dff\in\dff \nnn(V)$\dss is anisotropic and hence is non-degenerate.\oss
Indeed,\oss if\qss $n\fff(v)\qff \neq\qff 0$\dnsp,\oss
then every non-zero element\dss $w$\dss of\dss $V$\dss is equal to\dss $z\dff v$\dss for some non-zero\qss
$z\dff\in\dff \kkk$\nnsp.\oss
Since\dss $N$\dss is anisotropic,\oss
together with\qss (\ref{n-forms})\qss this implies that\qss $n\fff(w)\qff \neq\qff 0$\dnsp.\oss

The polarization of any norm-like form\dss $n$\dss is closely related to 
the polarization\qss
$\langle\fff\halfff z\fff,\qff u\fff\rangle$\qss
of\dss $N$\nnsp.\oss
Let\qss $v\dff\in\dff V$\qss
and\pss $z\fff,\pff u\dff\in\dff \kkk$\nnsp.\qff\oss
Then\oss
\begin{equation}
\label{n-form-polar}
\quad
(\fff z\dff v\fff,\qff u\dff v\fff)_n
\off =\off
\langle\fff\halfff z\fff,\qff u\fff\rangle\qff n\fff(v)\dff.
\end{equation}
Indeed,\oss
the maps\qss
$z\dff \longmapsto\dff n\fff(z\dff v)$\qss
and\qss
$z\dff \longmapsto\dff N\fff(z)\trf n\fff(v)$\qss
are quadratic forms on $\kkk$\nnsp.\oss
By\qss (\ref{n-forms})\qss they are equal.\oss
Hence their polarizations are equal also.\oss
But the left and the right hand sides of\qss (\ref{n-form-polar})\qss
are nothing else but these polarizations.\oss

\mypar{Lemma.}{zzbar-n-forms}
\emph{A quadratic form\dss $n$\trs on\qss $V$\qss is norm-like if and only if
\begin{equation}
\label{skew-self-adjoint}
\quad
(\fff z\dff v\fff,\qff w\fff)_n
\off =\off
(\fff v\fff,\qff  \oo{z}\qff w\fff)_n
\end{equation}
for every\qss $v\fff,\pff w\dff\in\dff V$\qss and every\qss
$z\dff\in\dff \kkk$\nnsp.\oss}

\proof\qss 
Suppose that\dss $n$\dss is norm-like.\oss
We may assume that\qss $v\qff \neq\qff 0$\dnsp.\oss
Then $v$ is a basis of $V$ over $\kkk$ and hence\qss
$w\qff =\qff u\dff v$\qss for some\qss $u\dff\in\dff \kkk$\nnsp.\oss
The identity\qss (\ref{n-form-polar})\qss implies that\qss
\[
\quad
(\fff z\dff v\fff,\qff w\fff)_n
\off =\off
(\fff z\dff v\fff,\qff u\dff v\fff)_n
\off =\off
\langle\fff\halfff z\fff,\qff u\fff\rangle\qff n\fff(v)
\hspace*{1.5em}\mbox{ and }\hspace*{1.5em}
\]

\vspace*{-3\bigskipamount}
\[
\quad
(\fff v\fff,\qff  \oo{z}\qff w\fff)_n
\off =\off
(\fff v\fff,\qff  \oo{z}\qff u\dff v\fff)_n
\off =\off
\langle\fff 1\fff,\qff \oo{z}\qff u\fff\rangle\qff n\fff(v)\dff.
\]
On the other hand,\oss
$\dis
\tr\dff(z\qff\oo{u}\fff)
\off =\off
\tr\dff(\fff\oo{z}\qff u)
\off =\off
\tr\dff(1\cdot \oo{z}\qff u)$\oss
and hence\oss
\[
\quad
\langle\fff\halfff z\fff,\qff u\fff\rangle
\off =\off
\langle\fff 1\fff,\qff \oo{z}\qff u\fff\rangle\dff.
\]
By combining these equalities 
we see that\qss (\ref{skew-self-adjoint})\qss holds.\oss
Conversely\halfff,\pss if\qss (\ref{skew-self-adjoint})\qss holds,\oss
then
\[
\quad
(\fff z\dff v\fff,\pff z\dff v\fff)_n
\off =\off
(\fff v\fff,\pff \oo{z}\dff z\dff v\fff)_n
\off =\off
z\dff\oo{z}\qff (\fff v\fff,\pff v\fff)_n
\]
and hence\oss
$\dis
n\fff(z\dff v)
\off =\off
N\fff(z)\dff n\fff(v)$\oss
for all\qss $z\dff\in\dff \kkk$\nnsp,\qss $v\dff\in\dff V$\dnsp.\oss  \eproof

\myuppar{The anti-norm-like quadratic forms.}
A form\qss $q\dff\in\dff Q\fff(V)$\qss is said to be\qss \emph{anti-norm-like}\qss if
\begin{equation}
\label{a-forms}
\quad
q\fff(\rho\dff v)
\off =\off
-\qff N\fff(\rho)\trf q\fff(v)
\off =\off
\rho^{\fff 2}\dff q\fff(v)
\end{equation}
for all\qss $v\dff\in\dff V$\qss and\qss $\rho\dff\in\dff \bkk$\nnsp.\oss
Let\dss $\aaa\fff(V)$\dss be the space of all anti-norm-like forms on\dss $V$\dnsp.\oss

As in the case of norm-like forms,\oss
the identity\qss (\ref{a-forms})\qss can be extended to the polarizations.\oss
Let\dss $q$\dss is an anti-norm-like form.\oss
Let\qss $v\fff,\pff w\dff\in\dff V$\qss and let\qss
$\rho\dff\in\dff \bkk$\nnsp.\oss
Then
\begin{equation}
\label{a-form-polar}
\quad
(\fff \rho\dff v\fff,\pff \rho\dff w\fff)_q
\off =\off 
-\qff N\fff(\rho)\dff 
(\fff v\fff,\pff w\fff)_q
\off =\off
\rho^{\fff 2}\dff
(\fff v\fff,\pff w\fff)_q
\end{equation}
Indeed,\oss the maps\qss
$v\dff \longmapsto\dff q\fff(\rho\dff v)$\qss
and\qss
$v\dff \longmapsto\dff -\qff N\fff(\rho)\dff q\fff(v)$\qss
are quadratic forms on $V$\dnsp.\oss
By\qss (\ref{a-forms})\qss they are equal.\oss
Hence their polarizations are equal also.\oss
But the left and the right hand sides of\qss (\ref{a-form-polar})\qss
are nothing else but these polarizations.\oss

\mypar{Lemma.}{mult-self-adjoint}
\emph{A quadratic form\dss $q$\trs on\qss $V$\qss is anti-norm-like if and only if
\begin{equation}
\label{self-adjoint}
\quad
(\fff z\dff v\fff,\pff w\fff)_q
\off =\off 
(\fff v\fff,\pff z\dff  w\fff)_q
\end{equation}
for every\qss $v\fff,\pff w\dff\in\dff V$\qss and every\qss
$z\dff\in\dff \kkk$\nnsp.\oss}

\proof\qss
Suppose that\dss $q$\dss is anti-norm-like.\oss
Every\qss $z\dff\in\dff \kkk$\qss
has the form\qss $z\qff =\qff a\qff +\qff \rho$\nnsp,\oss 
where\qss $a\dff\in\dff \kk$\qss and\qss $\rho\dff\in\dff \bkk$\nnsp.\off\oss 
If\oss $\rho\qff =\qff 0$\dnsp,\off\oss
then\oss $z\qff =\qff a\dff\in\dff \kk$\oss and hence\qss (\ref{self-adjoint})\qss holds.\oss
If\qss $\rho\qff \neq\qff 0$\dnsp,\oss
then\dss $\rho^{\fff -\dff 1}$\dss is defined
and\qss (\ref{a-form-polar})\qss
implies that\vspace*{\smallskipamount}
\[
\quad
(\fff \rho\dff v\fff,\pff w\fff)_q
\off =\off 
\rho^{\fff 2}\dff(\fff v\fff,\pff \rho^{\fff -\dff 1}\dff w\fff)_q 
\off =\off
(\fff v\fff,\pff \rho^{\fff 2}\dff\rho^{\fff -\dff 1}\dff w\fff)_q
\off =\off
(\fff v\fff,\pff \rho\dff w\fff)_q\qff.
\]

\vspace*{-0.75\bigskipamount}
On the other hand\halfff,\oss 
$\dis
(\fff a\dff v\fff,\pff w\fff)_q
\off =\off
(\fff v\fff,\pff a\dff w\fff)_q$\oss
because\qss $a\dff\in\dff \kk$\nnsp.\oss
Therefore\vspace*{\smallskipamount}
\[
\quad
(\fff a\dff v\fff,\pff w\fff)_q
\qff +\qff
(\fff \rho\dff v\fff,\pff w\fff)_q
\off =\off
(\fff v\fff,\pff a\dff w\fff)_q
\qff +\qff
(\fff v\fff,\pff \rho\dff w\fff)_q
\]

\vspace*{-0.75\bigskipamount}
It follows that\qss
$(\fff z\dff v\fff,\pff w\fff)_q
\qff =\qff 
(\fff v\fff,\pff z\dff  w\fff)_q$\nnsp,\qff\oss
i.e.\qss the identity\qss (\ref{self-adjoint})\qss holds.

Suppose now that\qss (\ref{self-adjoint})\oss holds.\oss
If\qss $v\dff\in\dff V$\qss and\qss $\rho\dff\in\dff \bkk$\nnsp,\oss
then\vspace*{\smallskipamount}
\[
\quad
(\fff \rho\dff v\fff,\pff \rho\dff v\fff)_q
\off =\off
(\fff \rho^{\fff 2}\dff v\fff,\pff v\fff)_q
\off =\off
\rho^{\fff 2}\dff (\fff v\fff,\pff v\fff)_q
\]

\vspace*{-0.75\bigskipamount}
and hence\qss $q\fff(\rho\dff v)\qff =\qff \rho^{\fff 2}\dff q\fff(v)$\dnsp.\qff\oss
It follows that\dss
$q$\dss is anti-norm-like.\oss \eproof

\mypar{Theorem.}{description-of-complement}
\emph{\dnsp$\aaa\fff(V)$\qss is equal to the orthogonal complement of\pss $\nnn\fff(V)$\qss
in\qss $Q\fff(V)$\dnsp.\oss}

\proof\qss Suppose that\qss $n\dff\in\dff \nnn\fff(V)$\dnsp,\oss $n\qff \neq\qff 0$\dnsp.\oss
Let\qss $v\dff\in\dff V$\dnsp,\qss $v\qss \neq\qss 0$\dnsp.\oss
Since $n$ is non-zero,\oss the property\qss (\ref{n-forms})\qss
implies that\qss $n\fff(v)\qff \neq\qff 0$\dnsp.\oss
Let\qss $\rho\dff\in\dff \bkk$\nnsp.\oss
Then\qss (\ref{n-form-polar})\qss implies that
\[
\quad
(\fff \rho\dff v\fff,\qff v\fff)_n
\off =\off
\langle\fff\halfff \rho\fff,\qff 1\fff\rangle\qff n\fff(v)
\off =\off
\tr\fff(\rho)\qff n\fff(v)
\off =\off
0
\]
and hence\qss $\rho\dff v$\qss 
is orthogonal to\qss $v$\qss with respect to $n$\nnsp.\oss

By Theorem\qss \ref{det-orthogonality}\qss 
the forms\qss $n\fff,\pff q$\qss are orthogonal if and only if
\[
\quad
n\fff(v)\dff q\fff(\rho\dff v)\qff +\qff q\fff(v)\trf n\fff(\rho\dff v)
\off =\off 0\dff.
\]
By the property\qss (\ref{n-forms})\qss
this equality is equivalent to
\[
\quad
n\fff(v)\dff q\fff(\rho\dff v)\qff +\qff q\fff(v)\trf N\fff(\rho)\dff n\fff(v)
\off =\off 0\dff.
\]
Since\qss $n\fff(v)\qff \neq\qff 0$\dnsp,\oss
the last equality is equivalent to\oss
\[
\quad
q\fff(\rho\dff v)\qff +\qff q\fff(v)\trf N\fff(\rho)
\off =\qff 0
\]
and hence to\oss
$\dis
q\fff(\rho\dff v)
\off =\off
-\qff N\fff(\rho)\dff q\fff(v)$\dnsp.\oss

It follows that\dss $n\fff,\pff q$\dss are orthogonal 
if and only if\dss $q$\dss is anti-norm-like.\oss  \eproof

\mypar{Theorem.}{direct-sum-decomposition} 
$\dis
Q\fff(V)\off =\off \nnn(V)\qff \oplus\qff \aaa\fff(V)$\dnsp.\oss

\proof\qss
Suppose that\oss 
$\dis
q\dff\in\dff \nnn\fff(V)\qff\cap\qff \aaa\fff(V)$\dnsp.\oss
Let\qss $\rho\dff\in\dff \bkk$\nnsp,\qss $\rho\qff \neq\qff 0$\dnsp.\oss
Then\oss 
\[
\quad
q\fff(\rho\dff v)\off =\off N\fff(\rho)\dff q\fff(v)
\hspace*{1.5em}\mbox{ and }\hspace*{1.5em}  
q\fff(\rho\dff v)\off =\off -\qff N\fff(\rho)\dff q\fff(v)
\]
and hence\qss
$\dis
N\fff(\rho)\dff
q\fff(v)\qff =\qff 0$\qss 
for all\qss $v\dff\in\dff V$\dnsp.\oss
But if\qss $\rho\qff \neq\qff 0$\dnsp,\oss
then\qss $N\fff(\rho)\qff \neq\qff 0$\dnsp.\oss
It follows that\qss $q\fff(v)\qff =\qff 0$\qss for all\dss $v$\nnsp,\oss
i.e.\qss $q\qff =\qff 0$\dnsp.\oss
Therefore\oss
\[
\quad
\nnn\fff(V)\qff\cap\qff \aaa\fff(V)
\off =\off
0
\]
In view of Theorem\qss \ref{description-of-complement},\oss
this implies\oss
$\dis
Q\fff(V)\off =\off N\fff(V)\qff\oplus\qff \aaa\fff(V)$\dnsp.\oss  \eproof

\mysection{The\qss structure\qss of\qss a\qss $\kkk$\dnsp-vector\qss space\qss 
on\qss $\aaa\fff(V)$}{a-multiplication}

\vspace*{\medskipamount}
\myuppar{The multiplication maps.}
For\qss $z\dff\in\dff \kkk$\qss
let\oss 
$\dis
m_{\dff z}\dff\colon\dff V\toto V$\oss 
be the map
\[
\quad
m_{\dff z}\dff \colon\dff v\qff \longmapsto\qff z\dff v\dff.
\]
We consider\dss $m_{\dff z}$\dss as a map linear over\dss $\kk$\dss
(although it is linear over\dss $\kkk$\dss also).\oss
Then\oss 
\begin{equation}
\label{determinant-of-multiplication}
\quad
\det\dff m_{\dff z}
\off =\off
N\fff(z)\dff.
\end{equation}
If\qss $V\qff =\qff \kkk$\nnsp,\oss
this is a well known result of the Galois theory\halfff.\oss
Since $V$ is iso\-morphic to $\kkk$ as a vector space over $\kkk$\nnsp,\oss
this special case implies the general one.\oss 

\myuppar{The anti-norm-like quadratic forms and the multiplication maps.}
Suppose that\dss $q$\dss is an anti-norm-like quadratic form on\dss $V$\dss 
and let\qss $z\dff\in\dff \kkk$\nnsp.\oss
Lemma\qss \ref{mult-self-adjoint}\qss implies that\dss $m_{\dff z}$\dss
is symmetric with respect to\dss $q$\nnsp,\oss and hence the product\dss
$m_{\dff z}\nsp\cdot\nsp q$\dss is defined.\oss

\mypar{Lemma.}{mult-z-a-form}
\emph{If\dss $q$\dss is an anti-norm-like quadratic form on\dss $V$\dss 
and\qss $z\dff\in\dff \kkk$\nnsp,\oss
then the product\qss $m_{\dff z}\nsp\cdot\nsp q$\qss 
is also anti-norm-like.}

\proof\qss Let\qss $\rho\dff\in\dff \bkk$\nnsp.\oss
By applying\qss (\ref{a-form-polar})\qss to\qss $z\dff v\fff,\pff v$\qss
in the roles of\qss $v\fff,\pff w$\nnsp,\oss we see that
\[
\quad
(m_{\dff z}\nsp\cdot\nsp q)\dff(\rho\dff v)
\off =\off
(\fff z\dff (\rho\dff v)\fff,\pff \rho\dff v\fff)_q
\off =\off
(\fff \rho\dff (z\dff v)\fff,\pff \rho\dff v\fff)_q
\]

\vspace*{-3\bigskipamount}
\[
\quad
\phantom{(m_{\dff z}\nsp\cdot\nsp q)\dff(\rho\dff v)
\off =\off
(\fff (z\dff \rho\dff v)\fff,\pff \rho\dff v\fff)_q
\off }
=\off
-\qff N\fff(\rho)\dff
(\fff z\dff v\fff,\pff v\fff)_q
\off =\off
-\qff N\fff(\rho)\qff
(m_{\dff z}\nsp\cdot\nsp q)\dff(v)\dff.
\]
It follows that\dss $m_{\dff z}\nsp\cdot\nsp q$\dss is anti-norm-like.\oss  \eproof

\mypar{Lemma.}{det-mult}
\emph{If\dss $q$\dss is an anti-norm-like quadratic form on\dss $V$\dss 
and\qss $z\dff\in\dff \kkk$\nnsp,\oss then 
\begin{equation*}
\quad
\det\qff (m_{\dff z}\nsp\cdot\nsp q)
\off =\off
N\fff(z)\dff \det q
\end{equation*}
for every determinant map\qss $\det\dff\colon\dff Q\fff(V)\toto \kk$\nnsp.\oss}

\proof\qss 
It is sufficient to combine\qss (\ref{det-dot-quadratic})\qss with\qss
(\ref{determinant-of-multiplication}).\oss  \eproof

\myuppar{The structure of a vector space over\dss $\kkk$\dss
on\dss $\aaa\fff(V)$\dnsp.}
Let us define the multiplication\qss 
$(z\fff,\pff q)\dff \longmapsto\dss z\fff q$\qss of forms\qss $q\dff\in\dff \aaa\fff(V)$\qss
by elements\qss $z\dff\in\dff \kkk$\qss by the formula
\[
\quad
z\fff q
\off =\off 
m_{\qff \oo{z}}\cdot\nsp q\dff.
\]
This turns $\aaa\fff(V)$ into a vector space over\dss $\kkk$\nnsp,\oss
as a trivial verification shows.\oss
By Lemma\qss \ref{det-mult}\qss any determinant map defines a norm-like 
quadratic form on\qss $\aaa\fff(V)$\dnsp.\oss

The multiplication\dss $m_{\qff \oo{z}}$\dss is used instead of\dss
the more natural\dss $m_{\dff z}$\dss 
in order to avoid\dss $\oo{z}$\dss in Lemma\qss \ref{homv-av-linearity}.\oss
This makes the identification maps from Section\qss \ref{comparing}\qss
linear over\dss $\kkk$\nnsp.\oss

\mysection{Conformal\qss structures}{conformal-structures}

\vspace*{\medskipamount}
\myuppar{Conformal structures.}
A\dss \emph{conformal structure}\dss on a vector space $U$ over $\kk$ is defined
as a non-zero quadratic form on $U$
considered up to multiplication by a non-zero element of $\kk$\nnsp.\oss
The conformal structure determined by\qss
$q\dff\in\dff Q\fff(U)$\qss is called the\dss \emph{conformal class}\qss of $q$
and is denoted by $[\fff q\dff]$\nnsp.\oss
It is called\qss
\emph{non-degenerate}\qss if the form $q$ is non-degenerate.\oss
The set of all conformal structures on $U$ is nothing else but the\qss
\emph{projective space}\dss $\mathbb{P}\hnsp Q\fff(U)$ associated with the\qss 
vector space\dss $Q\fff(U)$\dss of\ quadratic forms on $V $\dnsp.\oss

\myuppar{Conformal structures in dimension $1$ over $\kkk$\nnsp.}
The conformal class\qss $c_{\fff V}\qff =\qff [\fff n\fff]$\qss 
of a non-zero\dss $n\dff\in\dff \nnn(V)$\dss
does not depend on the choice of\dss $n$\dss and is called the\qss 
\emph{canonical conformal structure}\qss on $V$\dnsp.\oss 
A conformal structure on $V$ is called\qss 
\emph{exceptional}\qss if it is equal to the conformal class\dss $[\fff q\fff]$\dss
of a non-zero\dss $q\dff\in\dff \aaa\fff(V)$\nnsp,\oss
and\qss \emph{regular}\qss regular otherwise.\oss

Let\pss
$\dis
\mathbb{P}\hnsp \aaa\fff(V)\qff \subset\qff \mathbb{P}\hnsp Q\fff(V)$\pss
be the set of the exceptional conformal structures,\oss
i.e.\qss 
be the set of conformal classes\dss $[\fff q\dff ]$\dss
of quadratic forms\dss $q\dff\in\dff \aaa\fff(V)$\dnsp.\oss 
Its complement
\[
\quad
\mathbb{M}\fff(V)
\off =\off
\mathbb{P}\hnsp Q\fff(V)\qff \smallsetminus\qff \mathbb{P}\hnsp \aaa\fff(V)
\]
is the set of all regular conformal structures on $V$\dnsp.\oss
The set\qss $\mathbb{P}\hnsp \aaa\fff(V)$\qss is a line in
the projective plane\qss $\mathbb{P}\hnsp Q\fff(V)$\dnsp,\oss
and its complement\qss $\mathbb{M}\fff(V)$\qss is an affine plane.\oss

\myuppar{The structure of a vector space 
over\dss $\kkk$\dss on\qss $\mathbb{M}\fff(V)$\dnsp.}
Let\dss $n\dff\in\dff \nnn(V)$\dnsp,\qss $n\qff \neq\qff 0$\dnsp.\oss 
The map\qss 
$\dis
q\qff \longmapsto\qff [\fff n\qff +\qff q\fff]$\qss
is a bijection\qss
$\dis
\aaa\fff(V)\ttoo \mathbb{M}\fff(V)
$\dnsp.\oss 
One can use it to transfer the structure of a vector space over\dss $\kkk$\dss
from\dss $\aaa\fff(V)$\dss to\qss $\mathbb{M}\fff(V)$\dnsp.\oss
The resulting multiplication by elements\qss $z\dff\in\dff \kkk$\qss is given by the formula
\[
\quad
z\dff [\fff n\qff +\qff q\fff]
\off =\off
[\fff n\qff +\qff z\dff q\fff]\qff.
\]
While the bijection\qss
$\dis
\aaa\fff(V)\ttoo \mathbb{M}\fff(V)
$\qss
depends on the choice of\dss $n$\nnsp,\oss
an immediate verification shows that the resulting multiplication does not\halfff.\oss
Therefore,\oss
the above construction defines a canonical 
structure of a vector space 
over\dss $\kkk$\dss on\qss $\mathbb{M}\fff(V)$\dnsp.\oss
The canonical conformal structure\dss $c_{\fff V}$\dss on $V$ serves as the zero of this vector space.\oss

\myuppar{A canonical quadratic form on\qss $\mathbb{M}\fff(V)$\dnsp.}
Every point of\qss $\mathbb{M}\fff(V)$\qss has the form\qss $[\fff n\qff +\qff q\fff]$\nnsp,\oss
where\qss $n\dff\in\dff \nnn(V)$\dnsp,\qss $n\qff \neq\qff 0$\dnsp,\oss
and\qss $q\dff\in\dff \aaa\fff(V)$\dnsp.\qff\oss
Let us define a map\qss $\mathbb{D}\dff\colon\dff\mathbb{M}\fff(V)\toto \kk$\qss by
\begin{equation}
\label{mathbb-d}
\quad
\mathbb{D}\dff\colon\dff
[\fff n\qff +\qff q\fff]
\qff \longmapsto\qss
-\qff
\det\dff q\dff/\det\dff n\dff.
\end{equation}
If\oss 
$[\fff n'\qff +\qff q'\fff]
\off =\off
[\fff n\qff +\qff q\fff]$\nnsp,\oss
then\qss $n'\qff =\qff a\fff n$\nnsp,\qss 
$q'\qff =\qff a\fff q$\qss for some\qss $a\dff\in\dff \kk$\nnsp,\oss
and hence
\[
\quad
\det\dff q'\dff/\det\dff n'
\off =\off
\det\dff q\dff/\det\dff n\dff.
\]
It follows that\qss $\mathbb{D}$\qss is correctly defined.\oss
If we temporarily fix $n$ and identify\dss $\mathbb{M}\fff(V)$\dss
with\qss $\aaa\fff(V)$\qss by the map\qss 
$q\dff \longmapsto\dff [\fff n\qff +\qff q\fff]$\nnsp,\oss
then\qss $\mathbb{D}$\qss turns into the map\dss
$-\qff \det\nsp /\nsp\det\dff n$\dss 
restricted to\dss $\aaa\fff(V)$\dnsp.\oss
Since\dss $\det$\dss is a quadratic form,\oss
this implies that\dss $\mathbb{D}$\dss is a quadratic form on\dss $\mathbb{M}\fff(V)$\dnsp.\oss
Lemma\qss \ref{det-mult}\qss implies that\dss $\mathbb{D}$\dss
is norm-like.\oss

\mysection{Anti-linear\qss maps}{k-anti-linear}

\vspace*{\medskipamount}
\myuppar{Anti-linear maps.}
Let\qss $U\fff,\pff W$\qss be vector spaces over $\kkk$\dnsp,\oss
and let\dss $f$\dss be a map\qss $U\toto W$\qss linear over $\kk$\nnsp.\oss
The map\dss $f$\dss said to be\qss \emph{anti-linear}\qss over\dss $\kkk$\nnsp,\oss
or simply\qss \emph{anti-linear}\pss if
\begin{equation*}
\quad
f\fff(z\dff u)\off =\off \oo{z}\pff f\fff(u)
\end{equation*}
for all\qss $z\dff\in\dff \kkk$\dnsp,\qss $u\dff\in\dff U$\nnsp.\oss
Let\qss $\operatorname{Hom}_{\fff a}\fff(U\fff,\trf W)$\qss be
the space of all anti-linear maps\qss $U\toto W$\dnsp.

\mypar{Lemma.}{al-symmetric}
\emph{If\qss $f\dff\in\dff \operatorname{Hom}_{\fff a}\fff(V,\qff V)$\dss 
and\trs $n\fff\in\fff \nnn(V)$\dnsp,\oss
then\dss $f$\dss is symmetric with respect to $n$\nnsp.}

\proof\qss 
If\qss $v\fff,\pff w\dff\in\dff V$\qss and\qss $v\qff \neq\qff 0$\dnsp,\oss
then\qss $w\qff =\qff z\dff v$\qss for some\qss $z\dff\in\dff \kkk$\nnsp,\oss
and
\[
\quad
(\fff f\fff(v)\fff,\pff z\dff v\fff)_n
\off =\off
(\fff \oo{z}\qff f\fff(v)\fff,\pff v\fff)_n
\off =\off
(\fff f\fff(z\dff v)\fff,\pff v\fff)_n
\]
by Lemma\qss \ref{zzbar-n-forms}\qss and the anti-linearity of\qss $f$\nnsp.\qff\oss
It follows that\oss
\[
\quad
(\fff f\fff(v)\fff,\pff w\fff)_n
\off =\off
(\fff f\fff(w)\fff,\pff v\fff)_n
\off =\off
(\fff v\fff,\pff f\fff(w)\fff)_n
\]
and hence $f$ is symmetric with respect to\dss $n$\nnsp.\oss  \eproof

\mypar{Lemma.}{a-forms-maps}
\emph{Let\qss $n\dff\in\dff \nnn\fff(V)$\dnsp,\qss 
$n\qff \neq\qff 0$\dnsp.\qff\oss
If\pss $f\dff\colon\dff V\toto V$\qss is anti-linear\halfff,\oss
then\qss $f\nsp\cdot\nsp n$\qss is defined and\pss
$f\nsp\cdot\nsp n\dff\in\dff \aaa\fff(V)$\dnsp.\qff\oss
Conversely\halfff,\oss if\pss $q\dff\in\dff \aaa\fff(V)$\dnsp,\qff\oss 
then\pss $f\qff =\qff q/n$\qss
is anti-linear\halfff.\oss}

\proof\qss 
If\qss $f\dff\colon\dff V\toto V$\qss is anti-linear\halfff,\oss
then by Lemma\qss \ref{al-symmetric}\dss $f$\dss is symmetric with respect to\dss $n$\dss 
and hence\qss $f\nsp\cdot\nsp n$\qss is defined.\oss
Let\qss $q\qff =\qff f\nsp\cdot\nsp n$\nnsp.\oss
If\qss $z\dff\in\dff \kkk$\qss 
and\qss $v\fff,\pff w\dff\in\dff V$\dnsp,\oss 
then\oss
\[
\quad
(\dff \oo{z}\pff f\fff(v)\fff,\pff w\dff)_n
\off =\off
(\dff f\fff(z\dff v)\fff,\pff w\dff)_n
\off =\off 
(\fff z\dff v\fff,\pff w\fff)_q
\]
by anti-linearity of\dss $f$\dss and the definition of\dss $q$\nnsp,\oss
and\oss
\[
\quad
(\dff \oo{z}\pff f\fff(v)\fff,\pff w\dff)_n
\off =\off
(\dff f\fff(v)\fff,\pff z\dff w\dff)_n
\off =\off
(\fff v\fff,\pff z\dff w\fff)_q
\]
by Lemma\qss \ref{zzbar-n-forms}\qss
and the definition of\dss $q$\nnsp.\oss
It follows that\oss
$\dis
(\fff z\dff v\fff,\pff w\fff)_q
\off =\off
(\fff v\fff,\pff z\dff w\fff)_q$\oss
for all\qss $z\dff\in\dff \kkk$\qss and all\qss $v\fff,\pff w\dff\in\dff V$\dnsp.\oss
By Lemma\qss \ref{mult-self-adjoint}\qss
this implies that\qss $f\nsp\cdot\nsp n\qff =\qff q\dff\in\dff \aaa\fff(V)$\dnsp.\oss

Suppose now that\qss $q\dff\in\dff \aaa\fff(V)$\dnsp.\oss
Since $n$ is non-degenerate,\oss 
the map\oss
$\dis
f\qff =\qff 
q/n\dff\colon\dff V\toto V$\oss 
is defined and is symmetric with respect to\dss $n$\nnsp.\qff\oss
If\qss $z\dff\in\dff \kkk$\qss and\qss $v\fff,\pff w\dff\in\dff V$\dnsp,\qff\oss 
then\oss
\[
\quad
(\fff z\dff v\fff,\pff w\fff)_q
\off =\off
(\fff v\fff,\pff z\dff w\fff)_q
\]
by Lemma\qss \ref{mult-self-adjoint}\qss 
and hence\oss
$\dis
(\dff f\fff(z\dff v)\fff,\pff w\dff)_n
\off =\off
(\dff f\fff(v)\fff,\pff z\dff w\dff)_n
$\oss
by the definition of\dss $f$\nnsp.\oss
Since
\[
\quad
(\dff f\fff(v)\fff,\pff z\dff w\dff)_n
\qff =\qff
(\dff \oo{z}\qff f\fff(v)\fff,\pff w\dff)_n
\]
by Lemma\qss \ref{zzbar-n-forms},\oss
it follows that\oss
\[
\quad
(\dff f\fff(z\dff v)\fff,\pff w\dff)_n
\off =\off
(\dff \oo{z}\qff f\fff(v)\fff,\pff w\dff)_n
\]
for all\qss $z\dff\in\dff \kkk$\qss and all\qss $v\fff,\pff w\dff\in\dff V$\dnsp.\oss
Since\dss $n$\dss is non-degenerate and\dss $w$\dss is arbitrary\halfff,\oss
it follows that\qss $f\fff(z\dff v)\qff =\qff \oo{z}\qff f\fff(v)$\qss
for all\qss $z\dff\in\dff \kkk$\nnsp,\qss $v\dff\in\dff V$\dnsp,\qff\oss
i.e.\qss that\dss $f$\dss is anti-linear\halfff.\oss  \eproof

\mypar{Lemma.}{homv-av-linearity}
\emph{Let\qss $n\dff\in\dff \nnn\fff(V)$\dnsp,\qss 
$n\qff \neq\qff 0$\dnsp.\qff\oss
If\oss $f\dff\in\dff \homv$\oss and\oss $z\dff\in\dff \kkk$\nnsp,\qff\oss
then}\oss\vspace*{-1pt}
\[
\quad
(z\dff f)\nsp\cdot\nsp n
\off =\off
z\dff (f\nsp\cdot\nsp n)\dff.
\]

\vspace*{-\bigskipamount}
\proof\qss
Let\qss $q\qff =\qff f\nsp\cdot\nsp n$\qss
and\qss $p\qff =\qff (z\dff f)\nsp\cdot\nsp n$\dnsp.\oss
Let\qss $v\fff,\pff w\dff\in\dff V$\dnsp.\qff\oss
Then\oss
\[
\quad
(\fff v\fff,\pff w\fff)_{\fff z\fff q}
\off =\off
(\fff m_{\qff \oo{z}}\trf(v)\fff,\pff w\fff)_q
\off =\off
(\fff \oo{z}\qff v\fff,\pff w\fff)_q
\]
by the definition of the product\dss $z\fff q$\nnsp.\oss
On the other hand,\qss 
$(z\dff f)\fff(v)
\qff =\qff 
z\dff f\fff(v)
\qff =\qff
f\fff(\dff\oo{z}\qff v)
$\qss by the definition 
of the product\dss $z\fff f$\dss 
and the anti-linearity of\dss $f$\dss and hence
\[
\quad
(\fff v\fff,\pff w\fff)_{\fff p}
\off =\off
(\fff f\fff(\dff\oo{z}\qff v)\fff,\pff w\fff)_{\fff n}
\off =\off
(\fff \oo{z}\qff v\fff,\pff w\fff)_q
\]
by the definition of\trs $p$\dss and\dss $q$\nnsp.\qff\oss
It follows that\qss $p\qff =\qff z\fff q$\nnsp,\oss
i.e.\qss
$(z\dff f)\nsp\cdot\nsp n
\off =\off
z\dff (f\nsp\cdot\nsp n)$\dnsp.\oss  \eproof

\mysection{The\qss structure\qss of\pss $\homv$}{structure-of-a-hom}

\vspace*{\medskipamount}
\myuppar{The structure of a vector space over\dss $\kkk$\dss
on\dss $\homv$\dnsp.}
The product\qss 
$z\fff f$\qss of an element\qss $z\dff\in\dff \kkk$\qss with\qss $f\dff\in\dff \homv$\qss
is defined by the formula\vspace*{-1pt}
\[
\quad
z\fff f
\off =\off
m_{\dff z}\circ f
\dff.
\]

\vspace*{-1pt}
\vspace*{-\bigskipamount}
Obviously\halfff,\oss this operation turns $\homv$ into a vector space over\dss $\kkk$\nnsp.\oss
Since the dimension of\dss $V$\dss over\dss $\kkk$\dss is\dss $1$\nnsp,\oss
an anti-linear map\qss $V\toto V$\qss is determined by 
its value on any given non-zero\qss $v\dff\in\dff V$\dnsp.\oss
It follows that the dimension of\dss $\homv$ over\dss $\kkk$\dss is also\dss $1$\nnsp,\oss
and hence any non-zero\trs $g\dff\in\dff \homv$\trs is a basis of\dss $\homv$\dnsp.

The composition\qss $(f\fff,\pff g)\dff \longmapsto\dff f\circ g$\qss
is a\qss \emph{hermitian map}\pss\vspace*{-1pt}
\[
\quad
\homv\qff \times\qff \homv\ttoo\homv
\] 

\vspace*{-1pt}
\vspace*{-\bigskipamount}
in the sense that it is bilinear and\oss
$\dis
(z\fff f)\circ(u\fff g)
\off =\off
z\dff \oo{u}\pff (f\circ g)\dff
$\oss
for all\qss $z\fff,\pff u\dff\in\dff \kkk$\nnsp.\oss

\myuppar{A canonical bilinear form on\qss $\homv$\dnsp.}
Given\qss $f\fff,\pff g\dff\in\dff \homv$\dnsp,\qff\oss
let\vspace*{-1pt}
\[
\quad
\langle\fff f\fff,\pff g\fff\rangle
\off =\off
\Tr\dff (f\circ g)\bigl/2\qff,
\]

\vspace*{-1pt}
\vspace*{-\bigskipamount}
where the trace is taken over $\kk$\nnsp.\oss
Clearly\halfff,\qss 
$\langle\fff f\fff,\pff g\fff\rangle$\qss
is a bilinear form and since\vspace*{-1pt}
\[
\quad
\Tr\dff (f\circ g)
\off =\off
\Tr\dff (g\circ f)
\]

\vspace*{-1pt}
\vspace*{-\bigskipamount}
for all maps\qss 
$f\fff,\pff g\dff\colon\dff V\toto V$\qss
linear over $\kk$\nnsp,\oss 
this bilinear form is symmetric.\oss
Therefore\oss\vspace*{-1pt} 
\[
\quad
\mathcal{D}\dff\colon\dff f\qff\longmapsto\qff 
\langle\fff f\fff,\pff f\fff\rangle\bigl/2
\]

\vspace*{-1pt}
\vspace*{-\bigskipamount}
is a quadratic form on\qss $\homv$\qss 
having the form\qss
$\langle\fff f\fff,\pff g\fff\rangle$\qss 
as its polarization.\oss

\myuppar{Reflections.}
Let\dss $\sigma$\dss be the map $z\toto \oo{z}$\nnsp.\oss
Then\dss $\sigma$\dss is an anti-linear map\qss $\kkk\toto \kkk$\qss
and\qss $\sigma\circ\sigma\qff =\qff \id_{\dff \kkk\fff}$\dnsp.\oss
Since\qss $\sigma\fff(a)\qff =\qff a$\qss for\qss $a\dff\in\dff \kk$\qss
and\qss $\sigma\fff(\rho)\qff =\qff -\qff \rho$\qss for\qss $\rho\dff\in\dff \bkk$\nnsp,\oss
the determinant\dss $\det \sigma$\dss over\dss $\kk$\dss is equal to\dss $-\qff 1$\nnsp.\oss

Since\dss $V$\dss is isomorphic to\dss $\kkk$\dss as a vector space over\dss $\kkk$\nnsp,\oss
it follows that there exist anti-linear maps\qss $g\dff\colon\dff V\toto V$\qss such that\oss
$\dis
g\circ g\off =\off \id_{\dff V\fff}$\oss
and\oss
$\det\dff g\off =\off -\qff 1$\nnsp.\oss
We will call any such anti-linear map\qss $V\toto V$\qss a\qss \emph{reflection}.\oss

\mypar{Lemma.}{det-trace}
\emph{If\oss $f\dff\in\dff \homv$\dnsp,\qff\oss
then\oss
$\dis
f\circ f\qff
\off =\off
\mathcal{D}\dff(\dff f\dff)\dff \id_{\dff V\fff}$\oss and\oss
$\dis
\det\dff f
\off =\off
-\qff
\mathcal{D}\dff(\dff f\dff)$\dnsp,\oss
where the determinant is taken over\qss $\kk$\nnsp.\oss}

\proof\qss
Let\qss $g\dff\in\dff \homv$\qss be a reflection.\oss 
Then\qss $g\qff \neq\qff 0$\qss and hence\qss $g$\qss is a basis of\qss $\homv$\dnsp.\qff\oss
It follows that\oss $f\off =\off z\dff g$\oss
for some\qss $z\dff\in\dff \kkk$\qss 
and\dss hence
\[
\quad
f\circ f
\off =\off
(z\dff g)\circ (z\dff g)
\off =\off
z\dff\oo{z}\pff (g\circ g)
\off =\off
N\fff(z)\dff \id_{\dff V\fff}.
\]
Therefore\oss
$\mathcal{D}\dff(\dff f\dff)
\off =\off
\Tr\dff (\dff f\circ f\dff)/2
\off =\off
N\fff(z)$\dnsp.\qff\oss
It follows that\oss
$\dis
f\circ f\qff
\off =\off
\mathcal{D}\dff(\dff f\dff)\dff \id_{\dff V\fff}$\nnsp.\oss

On the other hand,\oss
$f\off =\off m_{\dff z}\circ g$\oss and\dss hence
Lemma\qss \ref{det-mult}\qss implies that
\[
\quad
\det\dff f
\off =\off
\left(\fff\det\dff m_{\dff z}\fff\right)\cdot \left(\fff\det\dff g\fff\right)
\off =\off
N\fff(z)\cdot(-\qff 1)
\off =\off
-\qff N\fff(z)\dff.
\]
Since we already proved that\oss
$\mathcal{D}\dff(\dff f\dff)
\off =\off
N\fff(z)$\dnsp,\oss
it follows that\oss
$\dis
\det\dff f
\off =\off
-\qff
\mathcal{D}\dff(\dff f\dff)$\dnsp.\oss  \eproof

\mysection{The\qss conformal\qss dilatation\qss and\qss Beltrami\qss forms}{conformal-and-beltrami}

\vspace*{\medskipamount}
\myuppar{Pull-backs of quadratic forms.}
Let $U\fff,\pff U'$ be
vector spaces over $\kk$ and $f\dff\colon\dff U\toto U'$ be a linear map.\qff\oss
Let\dss $q\dff \in\dff Q\fff(U')$\dnsp.\qff\oss
The\dss \emph{pull-back}\qss $f^{\dff *}\fff q\dff \in\dff Q\fff(U)$\dss 
is the quadratic form
\[
\quad
f^{\dff *}\fff q\qff\colon\qff u\qff \longmapsto\qff q\fff(f\fff(u))\dff.
\]
The pull-backs\dss $f^{\dff *}\fff q$\dss lead to a linear map\oss
$\dis
f^{\dff *}\qff\colon\qff
Q\fff(U')\ttoo Q\fff(U)$\dnsp.\oss
If $U''$ is one more vector space over $\kk$ and\qss 
$g\dff\colon\dff U'\toto U''$\qss is a linear map,\oss then\qss
$\dis
(g\circ\fff f)^*\qff =\qff f^{\dff *}\circ\qff g^{\fff *}$\dnsp.\oss
If $f$ is an isomorphism,\oss
then\dss $f^{\dff *}$\dss is an isomorphism also.\oss

If\qss $p\qff\in\dff Q\fff(U')$\qss 
and\qss $f^{\dff *}\dff p\qff \neq\qff 0$\dnsp,\oss
then the conformal structure $[\dff f^{\dff *}\dff p\dff]$ depends only on $[\dff p\dff]$\nnsp.\oss
It is called the\qss \emph{pull-back}\qss of\dss $p$\dss by\dss $f$\dss
and is denoted by\qss $f^{\dff *}\dff [\dff p\dff]$\nnsp.\oss

\myuppar{The conformal dilatation.}
Let $W$ be another vector space of dimension $1$ over $\kkk$
and let\qss $f\dff\colon\dff V\toto W$\qss be a non-zero map linear over $\kk$\nnsp.\oss
The pull-back 
\[
\quad
\mathbold{c}\fff(f)
\off =\off
f^{\dff *}\fff c_{\fff W}\dff\in\dff \mathbb{P}\hnsp Q\fff(V)
\] 
of the canonical conformal structure\dss $c_{\fff W}$\dss on\dss $W$\dss
is called the\qss \emph{conformal dilatation}\qss of\dss $f$\nnsp.\oss

The map\dss $f$\dss is called\qss \emph{exceptional}\qss if the conformal structure\dss
$\mathbold{c}\fff(f)$\dss is exceptional\halfff,\oss
and\qss \emph{regular}\qss otherwise.\oss
If\dss $f$\dss is regular\halfff,\oss
then we consider\dss
$\mathbold{c}\fff(f)$\dss
as an element of\qss $\mathbb{M}\fff(V)$\dnsp.\oss

If the field extension\dss $\kkk/\fff\halfff\kk$\dss resembles enough
the classical case\dss $\cc/\fff\rr$\nnsp,\oss
then all non-zero $f$ are regular\halfff.\oss
For example,\oss
by using Theorem\qss \ref{main-regular}\qss
it is not hard to prove that 
this is the case\dss if\dss $\kk$\dss is an ordered field
and\qss $\kkk\pff =\pff \kk\dff(\rho)$\qss with\qss $\rho^{\fff 2}\qff <\qff 0$\dnsp.\oss 

The map\dss $f$\dss is called\qss \emph{conformal}\pss if\dss
the pull-back by\dss $f$\dss of the canonical conformal structure on\dss $W$\dss
is equal to the canonical conformal structure on\dss $V$\dnsp,\oss i.e.\qss if
\[
\quad
f^{\dff *}\fff c_{\fff W}
\off =\off
c_{\fff V}\dff.
\] 
Since\dss $c_{\fff V}$\dss serves as the zero of\qss $\mathbb{M}\fff(V)$\dnsp,\oss
the map\dss 
$f$\dss is conformal if and only if\dss $f$\dss is regular
and\dss $\mathbold{c}\fff(f)\qff =\qff 0\dff\in\dff \mathbb{M}\fff(V)$\dnsp.\oss
In general,\oss the conformal dilatation of\dss $f$\dss is a natural measure
of the distortion of the canonical conformal structure by $f$\nnsp.\oss

\myuppar{Linear maps over\dss $\kk$\dss as sums of linear and anti-linear over\dss $\kkk$\dss maps.}
Every map\qss $f\dff\colon\dff V\toto W$\qss linear over\qss $\kk$\qss
admits a unique presentation as a sum
\[
\quad
f\off =\off Lf\qff +\qff Af
\]
of a map\dss $Lf$\dss linear over\dss $\kkk$\dss and a map\dss $Af$\dss
anti-linear over\dss $\kkk$\nnsp.\oss
Indeed,\oss
let\qss $\rho\dff\in\dff \bkk$\nnsp,\qss $\rho\qff \neq\qff 0$\dnsp.\qff\oss
If\qss 
$f\off =\off Lf\qff +\qff Af$\oss 
is\dss such a presentation of\dss $f$\nnsp,\oss
then\vspace*{\medskipamount}
\begin{equation*}
\quad
f\fff(v)
\off =\off Lf\fff(v)\qff +\qff Af\fff(v)\dff,
\hspace*{1.5em}
f\fff(\rho\dff v)
\off =\off
\rho\dff Lf\fff(v)\qff -\qff \rho\dff Af\fff(v)\dff,
\end{equation*}
and hence
\begin{equation*}
\quad
Lf\fff(v)\off =\off
\left(f\fff(v)\qff -\qff \rho^{\fff -\dff 1}\fff f\fff(\rho\dff v)\right)\dnsp\bigl/2\qff,
\hspace*{1.5em} 
Af\fff(v)\off =\off
\left(f\fff(v)\qff +\qff \rho^{\fff -\dff 1}\fff f\fff(\rho\dff v)\right)\dnsp\bigl/2
\end{equation*}

\vspace*{-\medskipamount}
for all\qss $v\dff\in\dff V$\dnsp.\oss
This proves uniqueness.\oss
Conversely\halfff,\oss one can define\qss $Lf$\qss and\qss $Af$\qss 
by the last displayed formulas.\oss
Since\dss $\bkk$\dss has dimension $1$ over $\kk$\nnsp,\oss
this definition does not depend on the choice of $\rho$\nnsp.\oss
We leave to the reader the verification that so defined maps\qss
$Lf$\qss and\qss $Af$\qss
are linear and anti-linear over $\kkk$ respectively\halfff.\oss

\myuppar{Beltrami forms.}
Let\qss $f\dff\colon\dff V\toto W$\qss is a map linear over $\kk$
such that\qss $Lf\qff \neq\qff 0$\dnsp.\oss
Then
\[
\quad
\mu_{\dff f}
\off =\off
(Lf)^{\fff -\dff 1} (Af)\qff \in\qff \homv
\]
is called the\qss \emph{Beltrami form}\qss of\qss $f$\nnsp.\oss
The map\dss $f$\dss is linear over\dss $\kkk$\dss if and only if\qss
$\mu_{\dff f}\qff =\qff 0$\dnsp.\oss
In general,\oss
the Beltrami form of\dss $f$\dss is a natural measure of deviation of\dss $f$\dss
from being linear over\dss $\kkk$\nnsp.\oss
It is a direct generalization of the classical
Beltrami forms from the theory of quasi-conformal mappings.\oss
See,\oss for example,\oss \cite{h},\oss Section\qss 4.8.\oss

\mysection{An\qss identification\qss of\oss $\mathbb{M}\fff(V)$\pss with\pss $\homv$}{identification}

\vspace*{\medskipamount}
\myuppar{An identification of\qss 
$\mathbb{M}\fff(V)$\qss with\qss $\homv$\dnsp.}
Every conformal structure in\qss $\mathbb{M}\fff(V)$\qss
has the form\qss $[\fff n\qff +\qff q\fff]$\qss
for some non-zero\qss $n\dff\in\dff \nnn(V)$\qss 
and some\qss $q\dff\in\dff \aaa\fff(V)$\dnsp.\qff\oss
Let  
\begin{equation}
\label{mv-to-homv}
\quad
\mathbb{M}\fff(V)\ttoo \homv
\end{equation}
be the map defined by the rule\oss
$\dis
[\fff n\qff +\qff q\fff]\qff \longmapsto\qff q/n$\nnsp.\oss
Since,\oss obviously\halfff,\oss
$\dis
a\fff q/a\fff n
\off =\off 
q/n$\oss 
for all non-zero\dss 
$a\dff\in\dff \kk$\nnsp,\qff\oss
this map is well defined.\oss

Let us choose some\qss 
$n\dff\in\dff \nnn(V)$\dnsp,\qss $n\qff \neq\qff 0$\dnsp.\qff\oss
Let  
\begin{equation}
\label{homv-to-mv}
\quad
\homv\ttoo \mathbb{M}\fff(V)
\end{equation}
be the map defined by the rule\oss
$\dis
f\qff \longmapsto\qff [\fff n\qff +\qff f\nsp\cdot\nsp n\fff]$\nnsp.\qff\oss
Since,\oss obviously\halfff,\oss
$\dis
f\nsp\cdot\nsp (a\fff n)\off =\off a\fff (f\nsp\cdot\nsp n)$\oss
and hence\oss
$\dis
[\fff a\fff n\qff +\qff f\nsp\cdot\nsp (a\fff n)\fff]
\off =\off
[\fff a\dff (n\qff +\qff f\nsp\cdot\nsp n)\fff]
\off =\off
[\fff n\qff +\qff f\nsp\cdot\nsp n\fff]$\oss for all non-zero\qss $a\dff\in\dff \kk$\nnsp,\oss
this map in fact does not depend on the choice of $n$\nnsp.\oss

Since the operations\qss $q\dff \longmapsto q/n$\qss and 
$f\dff \longmapsto f\nsp\cdot\nsp n$ are mutually inverse,\oss 
these two maps between\dss
$\mathbb{M}\fff(V)$\dss and\dss $\homv$\dss
are mutually inverse bijections.\oss
Lemma\qss \ref{homv-av-linearity}\qss implies that the second map
is linear over $\kkk$\nnsp.\oss
It follows that the first map is linear over $\kkk$ also,\oss
and hence both of these maps are isomorphisms of vector spaces over $\kkk$\nnsp.\oss
We will identify\qss
$\mathbb{M}\fff(V)$\qss and\qss $\homv$\qss
by these isomorphisms.\oss

\myuppar{The identification of the quadratic forms\dss 
$\mathbb{D}$\dss and\dss $\mathcal{D}$\dnsp.}
Let\qss $\det\dff\colon\dff Q\fff(V)\toto \kk$\qss be a determinant map.\oss
Let\qss 
$n\dff\in\dff \nnn(V)$\dnsp,\qss $n\qff \neq\qff 0$\dnsp.\oss
Suppose that\qss $f\dff\in\dff \homv$\qss
and let\qss $q\qff =\qff f\nsp\cdot\nsp n$\nnsp.\oss
Then\dss $f$\dss is identified with the conformal class\oss 
$[\fff n\qff +\qff f\nsp\cdot\nsp n\fff]
\off =\off
[\fff n\qff +\qff q\fff]$\nnsp.\oss

By applying\qss  (\ref{det-dot-quadratic})\qss to\dss $f$\dss and to\dss $n$\dss
in the role of\dss $q$\nnsp,\oss we see that\oss
$\det\dff q
\off =\off
\det\dff f \cdot \det\dff n$\nnsp.\oss
Since\oss
$\dis
\det\dff f
\off =\off
-\qff
\mathcal{D}\dff(\dff f\dff)$\oss
by Lemma\qss \ref{det-trace},\qff\oss
it follows that
\[
\quad
\mathcal{D}\dff(\dff f\dff)
\off =\off
-\qff 
\det\dff q\dff/\det\dff n\dff. 
\]
On the other hand,\oss
by the definition\qss (\ref{mathbb-d})\qss the quotient\qss
$-\qff \det\dff q\dff/\det\dff n$\qss is nothing else but
the value of\qss $\mathbb{D}$\qss on the conformal class\qss
$[\fff n\qff +\qff q\fff]$\qss identified with\dss $f$\nnsp.\oss

It follows that the pull-back of the quadratic form\qss $\mathcal{D}$\qss 
by the map\qss (\ref{homv-to-mv})\qss is equal to\qss $\mathbb{D}$\dnsp.\oss
Since the map\qss (\ref{mv-to-homv})\qss is equal to the inverse of\qss (\ref{homv-to-mv}),\oss
this implies that also the pull-back of the quadratic form\qss $\mathbb{D}$\qss 
by the map\qss (\ref{mv-to-homv})\qss is equal to\qss $\mathcal{D}$\dnsp.\oss

In other words,\oss
the identification of\qss $\mathbb{M}\fff(V)$\qss and\qss $\homv$\qss
identifies the canonical quadratic forms\qss $\mathbb{D}$\qss and\qss $\mathcal{D}$\qss
on these spaces.\oss

\mysection{Comparing\oss $\mathbold{c}\fff(f)$\oss and\oss $\mu_{\dff f}$}{comparing}

\vspace*{\medskipamount}
\myuppar{The framework.}
Let\qss $V\fff,\pff W$\qss be two vector spaces over\qss $\kkk$\nnsp,\oss
and\dss let\qss 
$f\dff\colon\dff V\toto W$\qss be a map linear over\qss $\kk$\nnsp.\oss 
Suppose that\qss $Lf\qff \neq\qff 0$\qss
and let\qss
$\mu
\qff =\qff 
\mu_{\dff f}$\qss 
be the Beltrami form of $f$\nnsp.

\mypar{Lemma.}{main-lemma}
\emph{The conformal dilatation\qss $\mathbold{c}\fff(f)$\qss of\qss $f$\qss
is equal to the conformal class of
\begin{equation}
\label{main-equation}
\quad
(\fff 1\qff +\qff \mathcal{D}\fff(\mu)\fff)\dff n
\pff +\pff
2\dff \mu\nsp\cdot\nsp n
\end{equation}
for any non-zero\pss
$n\dff\in\dff \nnn(V)$\nnsp.\oss}

\proof\qss
Since the conformal class of\qss 
$(\fff 1\qff +\qff \mathcal{D}\fff(\mu)\fff)\dff n
\pff +\pff
2\dff \mu\nsp\cdot\nsp n$\qss 
does not depend on the choice of\qss
$n\qff \neq\qff 0$\nnsp,\oss
it is sufficient to prove the claim
only for one particular\qss $n\qff \neq\qff 0$\dnsp.\oss

Let\qss $m\dff\in\dff \nnn(W)$\dnsp,\qss $m\qff \neq\qff 0$\dnsp.\oss
By the definition,\qss
$\mathbold{c}\fff(f)$\qss
is equal to the conformal class of\qss
$f^{\dff *}\dff m$\nnsp.\oss 
Let\qss $n\qff =\qff (\fff Lf\fff)^{\fff *}\fff(m)$\dnsp.\oss
Then\qss $n\dff\in\dff \nnn(V)$\qss because\qss $Lf$\qss is linear over\dss $\kkk$\dss
and\qss $n\qff \neq\qff 0$\qss because\qss
$Lf\qff \neq\qff 0$\qss and\dss $m$\dss is anisotropic.\qff\oss 
It follows that\oss
$\dis
(\fff Lf^{\qff -\dff 1}\fff)^{\fff *}\dff n\off =\off m$\oss and\dss hence\vspace*{3pt}
\[
\quad
f^{\dff *}\dff m
\off =\off
(\fff Lf\fff)^{\fff *}\fff m\qff +\qff (Af\fff)^{\fff *}\fff m
\]

\vspace*{-3\bigskipamount}
\[
\quad
\phantom{f^{\dff *}\dff m
\off }
=\off
(\fff Lf\fff)^{\fff *}\dff (\fff Lf^{\qff -\dff 1}\fff)^{\fff *}\dff n
\qff +\qff 
(Af\fff)^{\fff *}\fff (\fff Lf^{\qff -\dff 1}\fff)^{\fff *}\dff n
\]

\vspace*{-3\bigskipamount}
\[
\quad
\phantom{f^{\dff *}\dff m
\off }
=\off
n\qff +\qff (\fff Lf^{\qff -\dff 1}\circ Af\fff)^{\fff *}\dff n
\off 
=\off
n\qff +\qff \mu^{\fff *}\dff n
\off =\off
(\fff\id\qff +\qff \mu\fff)^{\fff *}\dff n\qff.\oss
\]

\vspace*{-\bigskipamount}\vspace*{3pt}
Let\qss $p\qff =\qff f^{\dff *}\dff m $\nsp.\qff\oss
By the above calculation\qss 
$p\qff =\qff (\fff\id\qff +\qff \mu\fff)^*\dff n$\qss 
and\dss hence\vspace*{\smallskipamount}
\[
\quad
(\fff v\fff,\pff w\fff)_{\fff p}
\off =\off
\left(\fff(\fff\id\qff +\qff \mu\fff)(v)
\fff,\pff
(\fff\id\qff +\qff \mu\fff)(w)\fff\right)_n
\]

\vspace*{-3\bigskipamount}\vspace*{-1pt}
\[
\quad
\phantom{(\fff v\fff,\pff w\fff)_{\fff p}
\off }
=\off
(\fff v\qff +\qff \mu\dff(v)
\fff,\pff
v\qff +\qff \mu\dff(w)\fff)_n
\]

\vspace*{-3\bigskipamount}\vspace*{-1pt}
\[
\quad
\phantom{(\fff v\fff,\pff w\fff)_{\fff p}
\off }
=\off
(\fff v\fff,\pff w)_n
\qff +\qff
(\fff \mu\dff(v)\fff,\pff \mu\dff(v)\fff)_n
\qff +\qff
(\fff v\fff,\pff \mu\fff(w)\fff)_n
\qff +\qff
(\fff \mu\dff(v)\fff,\pff w\fff)\dff.
\]

\vspace*{-0.75\bigskipamount}
By Lemma\qss \ref{al-symmetric}\qss the map\dss $\mu$\dss 
is symmetric with respect to\dss $n$\nnsp.\oss
It follows that
\[
\quad
(\fff v\fff,\pff w\fff)_{\fff p}
\off =\off
(\fff v\fff,\pff w\fff)_n
\qff +\qff
(\fff \mu\circ\mu\dff(v)\fff,\pff \fff w\fff)_n
\qff +\qff
2\dff (\fff \mu\dff(v)\fff,\pff w\fff)_n\dff.
\]
Lemma\qss \ref{det-trace}\qss implies that\oss
$\dis
\mu\circ \mu\qff
\off =\off
\mathcal{D}\dff(\dff \mu\dff)\dff \id_{\dff V\fff}$\oss
and\dss hence
\[
\quad
(\fff v\fff,\pff w\fff)_{\fff p}
\off =\off
(\fff v\fff,\pff w\fff)_n
\qff +\qff
\mathcal{D}\fff(\mu)\dff(\fff v\fff,\pff \fff w\fff)_n
\qff +\qff
2\dff (\fff \mu\dff(v)\fff,\pff w\fff)_n\dff.
\]
It follows that\oss
$\dis
f^{\dff *}\dff m
\off =\off
p
\off =\off
(\fff 1\qff +\qff \mathcal{D}\fff(\mu)\fff)\dff n
\pff +\pff
2\dff \mu\nsp\cdot\nsp n$\nnsp.\oss
Since\qss
$\mathbold{c}\fff(f)$\qss
is equal to the conformal class of\qss
$f^{\dff *}\dff m$\nnsp,\oss 
this proves the lemma.\oss  \eproof

\mypar{Theorem.}{main-regular}
\emph{If\oss
$\dis
1\qff +\qff \mathcal{D}\fff(\mu)\off \neq\off 0$\dnsp,\qff\oss
then\qss $f$\pss is regular\halfff,\oss
$\mathbold{c}\fff(f)\dff\in\qff \mathbb{M}\fff(V)$\dnsp,\qff\oss
and}\vspace*{\smallskipamount}
\[
\quad
\mathbold{c}\fff(f)
\off\qff =\off\qff
\frac{2\dff \mu}{\fff 1\qff +\qff \mathcal{D}\fff(\mu)\fff} 
\]

\vspace*{-0.75\bigskipamount}
\emph{after the identification of\pss $\mathbb{M}\fff(V)$\qss
with\qss $\homv$\dnsp.\oss}

\proof\qss
If\oss
$\dis
1\qff +\qff \mathcal{D}\fff(\mu)\off \neq\off 0$\dnsp,\oss
then the quadratic form\qss (\ref{main-equation})\qss
does not belong to\qss
$\aaa\fff(V)$\qss
and the conformal class of the form\qss (\ref{main-equation})\qss
is equal to the conformal class of\vspace{\smallskipamount}
\begin{equation}
\label{main-equation-regular}
\quad
n
\pff +\pff 
\frac{2\dff \mu}{1\qff +\qff \mathcal{D}\fff(\mu)}
\cdot n\dff.
\end{equation}

\vspace*{-0.75\bigskipamount}
In view of Lemma\qss \ref{main-lemma},\oss
this implies that\oss 
$\mathbold{c}\fff(f)\dff\in\qff \mathbb{M}\fff(V)$\oss
and\qss $\mathbold{c}\fff(f)$ is equal to the conformal class of\qss 
(\ref{main-equation-regular}).\oss
It remains to notice that conformal class is identified with\vspace{\smallskipamount}
\[
\quad
\frac{2\dff \mu}{1\qff +\qff \mathcal{D}\fff(\mu)}
\qff \in\qff
\homv\dff.
\hspace*{1.5em}\mbox{ \eproof}
\]

\vspace*{-0.75\bigskipamount}
\mypar{Theorem.}{main-exceptional}
\emph{If\oss
$\dis
1\qff +\qff \mathcal{D}\fff(\mu)\off =\off 0$\dnsp,\qff\oss
then\qss $f$\pss is exceptional\halfff,\oss
$\mathbold{c}\fff(f)\dff\in\qff \mathbb{P}\hnsp \aaa\fff(V)$\dnsp,\qff\oss
and}\oss
\[ 
\mathbold{c}\fff(f)
\off =\off
[\dff 2\dff \mu\nsp\cdot\nsp n\dff]
\] 
\emph{for any non-zero\qss $n\dff\in\dff \nnn(V)$\dnsp.\oss}

\proof\qss 
Since\qss $2\dff \mu\nsp\cdot\nsp n\dff\in\dff \aaa\fff(V)$\qss 
by Lemma\qss \ref{a-forms-maps},\oss
this follows from Lemma\qss \ref{main-lemma}.\oss  \eproof

\begin{flushright}

January\qss 21,\oss 2017
 
https\halfff:/\!/{\hnsp}nikolaivivanov.com

E-mail\halfff:\oss nikolai.v.ivanov{\fff}@{\fff}icloud.com

\end{flushright}


\begin{thebibliography}{X}



\bibitem[A]{at} M.\qss Atiyah,\oss 
The non-existent complex $6$\dnsp-sphere,\oss
arXiv:1610.09366v1,\oss 2016,\oss 7\qss pp.

\bibitem[H]{h} J.H.\qss Hubbard,\oss \emph{Teichmüller Theory and Applications to 
Geometry,\oss Topology,\oss and Dynamics,\oss Volume\qss I:\oss Teichmüller Theory},\oss 
Matrix Editions,\oss
2006,\oss xx,\qss 461\qss pp.

\bibitem[I]{i-d} N.V.\qss Ivanov,\oss The geometric meaning of the conformal dilatation,\oss 
arXiv:1701.06259,\oss 2017,\oss 20\qss pp.


\end{thebibliography}
\end{document}